\documentclass[10pt]{article}
\usepackage{amssymb}
\usepackage{latexsym}
\usepackage{amsmath}
\usepackage{graphics}
\newcommand{\PSbox}[3]{\mbox{\rule{0in}{#3}\includegraphics{#1}\hspace{#2}}} 
\newtheorem{Theorem}{Theorem}[section]
\newtheorem{Definition}{Definition}[section]

\newtheorem{Lemma}[Theorem]{Lemma}
\newtheorem{Proposition}[Theorem]{Proposition}

\newtheorem{Corollary}[Theorem]{Corollary}

\newtheorem{Claim}[Theorem]{Claim}
\newtheorem{Fact}{Fact}[section]
\newtheorem{Remark}{Remark}[section]
\numberwithin{equation}{section}

\newenvironment{Proof}{\par\noindent{\em Proof}. }{\hfill $\Box$\par\indent}
\def\sqr#1#2{{\vcenter{\vbox{\hrule  height.#2pt
       \hbox{\vrule width.#2pt height#1pt \kern#1pt \vrule width.#2pt}
        \hrule height.#2pt}}}}
   
\let\epsilon=\varepsilon
\let\Bbb=\mathbb
\let\phi=\varphi
\def\){ \right) }
\def\({ \left( }
\def\[{ \left[ }
\def\]{ \right] }
\def\<{ \langle }
\def\>{ \rangle }
\let\ljunk=\{
\let\rjunk=\}
\def\{{\left\ljunk}
\def\}{\right\rjunk}
\def\p{\partial}

\def\Ker{\mbox{Ker}}

\def\Riem{{\cal R}{\mathrm i}{\mathrm e}{\mathrm m}}

\def\sign{{\mathrm s}{\mathrm i}{\mathrm g}{\mathrm n}}

\def\Crit{{\mathrm C}{\mathrm r}{\mathrm i}{\mathrm t}}
\def\cyl{{\mathit c}{\mathit y}\ell}
\def\ext{{\mathrm e}{\mathrm x}{\mathrm t}}

\def\cp{{\mathrm c}}

\def\const{{\mathrm c}{\mathrm o}{\mathrm n}{\mathrm s}{\mathrm t}.}
\newcommand{\orb}{{\mathrm o}{\mathrm r}{\mathrm b}}

\newcommand{\supp}{{\mathbf S}{\mathrm u}{\mathrm p}{\mathrm p}}

\def\lra{\longrightarrow}
\def\L{\Bbb L}
\def\Vol{{\mathrm V}{\mathrm o}{\mathrm l}}

\def\Im{{\mathrm I}{\mathrm m}\ \!\!}
\def\Ell{\mathsf L}

\newcommand{\R}{{\mathbf R}}
\newcommand{\B}{{\mathbf B}}

\newcommand{\Z}{{\mathbf Z}}
\newcommand{\ind}{{\mathrm i}{\mathrm n}{\mathrm d} }

\newcommand{\Sp}{{\mathsf S}}
\newcommand{\C}{{\mathbf C}}

 \def\Slash#1{
  \begin{picture}(5,6)(0,0)
  \put(-.7,-1.2){\line(5,6)5}
  \end{picture}
  \kern-.8em#1}
\def\dirac{ \ {\Slash \partial}}

\def\P{\mbox{\bf P}}
\def\ccd{\begin{picture}(5,3)(0,0)
  \put(2.5,2){\circle*{1.5}}
  \end{picture}}
\newenvironment{lproof}[1]{\par\vspace{2mm} \noindent{\em #1.} }{\hfill $\Box$\par\vspace{2mm}}

\title{The Yamabe invariants of orbifolds and cylindrical manifolds, 
and $L^2$-harmonic spinors}
\author{Kazuo Akutagawa$\ {}^*$, Boris Botvinnik \thanks{Partially supported by the Grants-in-Aid for Scientific Research (C), Japan Society for the Promotion of Science, No. 14540072.} }
\date{ \ }

\begin{document} 
\maketitle
\markboth{The Yamabe invariants of orbifolds and cylindrical manifolds}
{K. Akutagawa, B. Botvinnik}
\vspace{-14mm}

\begin{abstract}
We study the Yamabe invariants of cylindrical manifolds and 
compact orbifolds with a finite number of singularities, 
by means of conformal geometry 
and the Atiyah-Patodi-Singer $L^2$-index theory. 
For an $n$-orbifold $M$ with singularities $\Sigma_{\Gamma} = 
\{(\check{p}_1, \Gamma_1), \ldots, (\check{p}_s, \Gamma_s)\}$ 
(where each group $\Gamma_j<O(n)$ is of finite order), 
we define and study the \emph{orbifold Yamabe invariant} $Y^{\orb}(M)$. 
We prove that $Y^{\orb}(M)$ coincides with the corresponding 
$h$-$\emph{cylindrical Yamabe invariant}$  
$Y^{h\textrm{-}\cyl}(M \setminus \{\check{p}_1, \ldots, \check{p}_s\})$ 
defined by the authors \cite{AB2}, 
where $h = h_{\Gamma_j}$ is the standard metric 
on the slice $S^{n-1}/\Gamma_j$ of each end with infinity $\check{p}_j$. 
Using this, we show that 
$Y^{\orb}(M)$ is bounded by $Y(S^n) /d$ from above, 
where $d=\max_j|\Gamma_j|^{\frac{2}{n}}$. 
For a cylindrical $4$-manifold $X$ with a general slice metric $h$ on the end, 
we also establish a method for estimating the $h$-cylindrical Yamabe invariant 
$Y^{h\textrm{-}\cyl}(X)$ from above, 
in terms of the geometry and topology of $X$. 
We conclude by an explicit estimate of $Y^{h\textrm{-}\cyl}(X)$ 
for particular cylindrical $4$-manifolds $X$, 
including that of $Y^{\orb}(M)$ for $4$-orbifolds $M$.
\end{abstract}

\section{Introduction}\label{int}
In this paper, we study the Yamabe invariant of compact orbifolds 
with a finite number of singularities (i.e., the orbifold Yamabe invariant)
and its relationship to the Yamabe invariant of cylindrical manifolds 
(i.e., the cylindrical Yamabe invariant). 
In the $4$-dimensional case, 
we also establish a method for estimating the cylindrical Yamabe invariant 
(including the orbifold Yamabe invariant) from above,
by means of conformal geometry 
and the Atiyah-Patodi-Singer $L^2$-index theory. 

\noindent
{\bf \ref{int}.1. The orbifold Yamabe invariant.} 
First, we define the Yamabe invariant of such orbifolds, 
that is, the orbifold Yamabe invariant,  
and then we study some fundamental properties 
on the orbifold Yamabe invariant.  
An orbifold $M$ under consideration here 
is a relatively compact smooth
manifold outside of a finite number of singular points 
$\check{p}_1,\ldots, \check{p}_s \in M$, 
and near each point $\check{p}_j$ 
it is locally homeomorphic to the orbit space $\R^n/\Gamma_j$, 
where $\Gamma_j ( < O(n) )$ is a finite group acting freely 
on $\R^n\setminus \{0\}$. 
We use the notation 
$\Sigma_{\Gamma} = \{(\check{p}_1,\Gamma_1),\ldots,(\check{p}_s,\Gamma_s)\}$ 
and $\Sigma = \{\check{p}_1, \ldots, \check{p}_s\}$, 
$\Gamma = \{\Gamma_1, \ldots, \Gamma_s\}$ 
for the singularities of $M$. 
Throught this paper, we assume that $\dim M = n \geq 3$. 

There are natural \emph{orbifold Riemannian metrics} 
compatible with the orbifold structure. 
We denote the space of such metrics by $\Riem^{\orb}(M)$. 
Let $R_g$, $d\sigma_g$ and $\Vol_g(M)$ denote respectively 
the scalar curvature, the volume form and the volume 
corresponding to $g\in \Riem^{\orb}(M)$. 
To start with, we note that the (normalized) 
Einstein-Hilbert functional
$$
I: \Riem^{\orb}(M) \to \R, \ \ \ \mbox{where} \ \ \ \ 
g\mapsto \frac{\int_M R_g d\sigma_g}{\Vol_g(M)^{\frac{n-2}{n}}},
$$
has exactly the same property as in the case of compact 
manifolds: the set of critical points $\Crit(I)$ coincides with the
space of \emph{orbifold Einstein metrics} (see Proposition \ref{EH}).

We also observe that the functional $I$ 
restricted on each \emph{orbifold conformal class}
$[g]_{\orb}\in {\mathcal C}^{\orb}(M)$ has the same basic
properties as in the case of compact manifolds (cf.~\cite{KW}), 
where ${\mathcal C}^{\orb}(M)$ denotes 
the space of orbifold conformal classes on $M$. 
This leads us to the definitions of the
\emph{orbifold Yamabe constant} $Y_{[g]_{\orb}}(M)$ and the
\emph{orbifold Yamabe invariant} $Y^{\orb}(M)$ respectively:
$$
Y_{[g]_{\orb}}(M):= \inf_{\tilde{g}\in [g]_{\orb}} I(\tilde{g}),
\ \ \ \
Y^{\orb}(M):=\sup_{C\in {\mathcal C}^{\orb}(M)} Y_C(M).
$$

\noindent
{\bf \ref{int}.2. 
The orbifold and the cylindrical Yamabe invariants.} 
Next, we study the relationship 
between the oribifold Yamabe invariant 
and the corresponding $h$-cylindrical Yamabe invariant. 
Consider, for simplicity, the case when a compact orbifold $M$ 
has only one singularity $\Sigma_{\Gamma} = \{(\check{p},\Gamma)\}$.  
We observe that, in the category of smooth manifolds, 
the open manifold $X:=M\setminus\{\check{p}\}$ could be considered as 
the underlying smooth manifold of a \emph{cylindrical manifold} $(X,\bar{g})$ 
with the cylindrical end 
$\left((S^{n-1}/\Gamma)\times[1,\infty), \bar{g}=h_{\Gamma}+dt^2\right)$, 
where $(S^{n-1}/\Gamma,h_{\Gamma})$ is 
the corresponding spherical space form 
equipped with the standard metric $h_{\Gamma}$ of constant curvature one. 
As an object in the category of Riemannian manifolds, 
the orbifold $M$ equipped with an orbifold metric $g$ 
is quite different from the cylindrical manifold $(X,\bar{g})$. 
Even from the viewpoint of conformal geometry, 
there is no ``cylindrical metric'' within the conformal class $[g]$, 
in general. 
However, the following holds (Theorem \ref{s2_Th1}). 

\vspace{1mm} 

\noindent
{\bf Theorem A.}  
{\sl Let $M$ be a compact orbifold with one singularity 
$\Sigma_{\Gamma} = \{(\check{p},\Gamma)\}$. 
Then $Y^{\orb}(M)= Y^{h_{\Gamma}\hbox{-}\cyl}(M\setminus \{\check{p}\})$. 
Here, $Y^{h_{\Gamma}\hbox{-}\cyl}(X)$ denotes 
the $h_{\Gamma}$-cylindrical Yamabe invariant of $X$ 
(see Section~2.6 or \cite[Section~2]{AB2} 
for the definition of $Y^{h_{\Gamma}\hbox{-}\cyl}(X)$).} 
\vspace{1mm} 

\noindent
Using this, 
we prove the following estimate of the orbifold Yamabe invariant 
(Corollary \ref{s2_Cor1}).
\vspace{1mm} 

\noindent
{\bf Theorem B.}  
{\sl Let $M$ be a compact orbifold of $\dim M=n\geq 3$ with singularities 
$\Sigma_{\Gamma} = 
\{(\check{p}_1,\Gamma_1), \ldots,(\check{p}_{s}, \Gamma_{s})\}$. 
Then} 
$$ 
Y^{\orb}(M) = Y^{h_{\Gamma}\hbox{-}\cyl}
(M\setminus \{\check{p}_1,\ldots,\check{p}_{s}\}) 
\leq \min_{1 \leq j \leq s}\frac{Y(S^n)}{|\Gamma_j|^{\frac{2}{n}}}, 
$$
{\sl where $h_{\Gamma} = h_{\Gamma_j}$ on the slice $S^{n-1}/\Gamma_j$ 
of each end $(S^{n-1}/\Gamma_j) \times [1, \infty)$, 
and $|\Gamma_j|$ denotes the order of $\Gamma_j$.} 
\vspace{1mm} 

\noindent
{\bf \ref{int}.3. Harmonic spinors on cylindrical $4$-manifolds.}  
Our main goal is to obtain 
an estimate of the $h$-cylindrical Yamabe invariant 
for a general slice metric $h$ 
and, via Theorem A, the orbifold Yamabe invariant, 
by means of the Atiyah-Patodi-Singer $L^2$-index theory 
for cylindrical $4$-manifolds. 
Indeed, in Section~3, 
we will give some estimates for the $h$-cylindrical Yamabe invariant 
from above, which can be regarded as 
a natural generalization of those for the Yamabe invariant 
obtained in Gursky-LeBrun~\cite{GL}. 

Let $(X,\bar{g})$ be a cylindrical $4$-manifold modeled by $(Z,h)$ 
for a metric $h$ on the slice $Z$. 
It means that on the end $Z\times [1,\infty)$, 
the metric $\bar{g}(z,t)=h(z)+dt^2$ is cylindrical 
with respect to a product coordinate system $(z,t)$. 
Then we use the notation that $\p_{\infty}\bar{g}=h$. 
Now, we assume that there exists a spin$^c$-structure 
on $(X,\bar{g})$ given by an element 
$$
a \in H^2(X;\Z) \cap \Im \left[H^2_c(X;\R) \to H^2(X;\R)\right]
$$ 
which is not a torsion class, 
where $H^2_c(X;\R)$ denotes the second cohomology with compact support. 
Let $\Ell\to X$ be the $\C$-line bundle with $c_1(\Ell)=a$ 
and with the restriction $\Ell|_{Z\times [1,\infty)}$ being ``cylindrical''. 
Similarly to the case of compact manifolds, 
there exists a unique \emph{$L^2$-harmonic form 
$\zeta\in {\mathcal H}_{\bar{g}}^2(X)$} 
representing the cohomology class $a$. 
Note that the form $\zeta$ does not have compact support 
unless $a=0\in H^2(X;\R)$. 
However, there exists a sequence $\{\zeta_j\}_{j=1}^{\infty}$ 
of closed $2$-forms with compact support such that $[\zeta_j] = a$ 
and $\zeta_j\to \zeta$ in an appropriate topology.  
We choose $U(1)$-connections $A_j$ (with compact support) on $\Ell\to X$ 
such that 
$\frac{\sqrt{-1}}{2\pi}F_{A_j}=\frac{\sqrt{-1}}{2\pi}dA_j = \zeta_j$. 
Then, for each $j$, the associated twisted Dirac operator 
$$
\dirac_{A_j} =
\dirac_{A_j}^{\ +}: \Gamma(\Sp^+_{\C}(a))\to \Gamma(\Sp^-_{\C}(a)) 
$$ 
is defined, 
where $\Sp^{\pm}_{\C}(a):= \Sp^{\pm}_{\C}\otimes \Ell^{1/2}$ 
stands for the plus/minus spin bundle associated to $\Ell$. 
Note that the $L^2$-index of $\dirac_{A_j}$ is independent of $j$. 
We denote by $b_2^-(X)$ the total dimension of the negative eigenspaces 
of the intersection form on $H^2_c(X; \R)$. 
The following is the main result of Section~\ref{s3} (Theorem \ref{s3_T1}). 

\noindent
{\bf Theorem C.} {\sl Let $(X,\bar{g})$ be a cylindrical $4$-manifold
modeled by $(Z,h)$ with $b_2^-(X)=0$. 
Assume that the $L^2$-index of $\dirac_{A_j}$ is positive. Then} 
$$
Y^{h\hbox{-}\cyl}(X)\leq 
\min\{4\pi\sqrt{2a^2},\ Y^{\cyl}_{[h+dt^2]}(Z\times \R)\}. 
$$ 

The proof of Theorem~C is subtle. 
We modify the argument due to Gursky-LeBrun \cite{GL}, 
involving the Bochner technique 
and a generalization of the \emph{modified scalar curvature} 
to the case of noncompact manifolds.  
One of the difficulties in the proof is the following: 
It turns out that a suitable conformal metric 
$\check{g}=u^2\cdot \bar{g}$ (where $u\in L^{1,2}_{\bar{g}}(X)$) 
for this argument never provides a complete metric on $X$ 
(see Lemmas \ref{s3_L2}, \ref{s3_L4} 
and the estimate (\ref{s3_eq6}) in Section~\ref{s3}). 
The main technical difficulty arises from this point. 
In order to overcome the difficulty, we have to estimate, for instance, 
the function $u^{-1}|du|_{\bar{g}}$ uniformly on the cylindrical end 
$Z\times [1,\infty)$ (see Lemma \ref{s3_L6}). 

We remark that 
the Hodge theorem and the Atiyah-Singer index theorem 
still hold for compact orbifolds (under some modifications) 
(cf.~\cite{Kawasaki}, \cite{Du}). 
Hence, if one considers only compact $4$-orbifolds 
and their orbifold Yamabe invariants, 
the arguments for the estimate 
are more direct with the aid of these theorems. 
However, we establish here an estimate 
of the cylindrical Yamabe invariant for cylindrical 4-manifolds 
modeled by general compact Riemannian $3$-manifolds. 
Moreover, the combination of conformal geometry and analysis 
on cylindrical manifolds is of independent value and 
leads to new interesting results (cf.~\cite{AB1}, \cite{AB2}, \cite{ABKS}). 
Therefore, we rather deal with conformal geometry and analysis 
on cylindrical manifolds 
by means of the $L^2$-Hodge theory 
and the Atiyah-Patodi-Singer $L^2$-index theory. 

\noindent
{\bf \ref{int}.4. Examples.} 
Here we present particular cylindrical $4$-manifolds $X_{\ell}$ 
defined as follows. 
Let $\Ell\to \C\P^1$ be the anti-canonical $\C$-line bundle over $ \C\P^1$, 
and $\Ell_{\ell}:=\Ell^{\otimes\ell}\to \C\P^1$ for $\ell\geq 1$.  
Then, let $X_{\ell}$ denote 
the total space of the bundle $\Ell_{\ell}\to \C\P^1$, 
which will be naturally the underlying smooth open manifold 
(with a connected tame end) of a cylindrical manifold. 
For any cylindrical metric $\overline{g}$ on $X_{\ell}$, 
the cylindrical end of $X_{\ell}$ is reprented by 
$((S^{3}/\Gamma_{\ell}) \times [1, \infty), h + dt^2)$ 
with a metric $h$ on $S^{3}/\Gamma_{\ell}$, 
where $\Gamma_{\ell}:=\Z/\ell\Z$. 
Note that the one-point compactification 
$M_{\ell} := X_{\ell} \cup \{\check{p}_{\infty}\}$ of $X_{\ell}$ 
has a natural orbifold structure with one singularity 
$\Sigma_{\Gamma} = \{(\check{p}_{\infty}, \Gamma_{\ell})\}$. 
The main result of Section \ref{s4} is the following estimate 
(Theorem \ref{s4_T1}): 
\vspace{1mm}

\noindent
{\bf Theorem D.} {\sl 
With the above understood, 
let $\Riem^+(S^3/\Gamma_{\ell})$ denote the space of metrics 
of positive scalar curvature on $S^3/\Gamma_{\ell}$. 
\begin{enumerate}
\item[{\bf (1)}] 
For $\ell\geq 1$ and any metric $h \in \Riem^+(S^{3}/\Gamma_{\ell})$ 
homotopic to $h_{\Gamma_{\ell}}$ in $\Riem^+(S^3/\Gamma_{\ell})$, 
$$
0 < Y^{h\hbox{-}\cyl}(X_{\ell}) \leq 
\min\{4\pi(\ell+2)\sqrt{\frac{2}{\ell}},\ 
Y^{\cyl}_{[h+dt^2]}((S^{3}/\Gamma_{\ell})\times \R)\}. 
$$
\item[{\bf (2)}] 
For any metric $h \in \Riem(S^{3}/\Gamma_{\ell})$ 
which is sufficiently $C^2$-close to $h_{\Gamma_{\ell}}$, 
$$
0< Y^{h\hbox{-}\cyl}(X_{\ell})\ \leq \ 
\{ 
\begin{array}{rll} 
12\sqrt{2}\pi = Y(\C\P^2)\quad\ \ & \mbox{if} \ \ \ \ell=1, 
\\ 
Y^{\cyl}_{[h+dt^2]}((S^{3}/\Gamma_{\ell})\times \R) 
& \mbox{if} \ \ \ \ell\geq 2. 
\end{array} 
\right. 
$$
In particular, 
$0 < Y^{\orb}(M_{\ell}) = 
Y^{h_{\Gamma_{\ell}}\hbox{-}\cyl}(X_{\ell}) \leq  
8\sqrt{\frac{6}{\ell}}\ \pi$ 
\ for $\ell \geq 2$. 
\end{enumerate} 
} 

\noindent
{\bf \ref{int}.5. The plan of the paper.} 
In Section \ref{s2}, 
we define and study the orbifold Yamabe constant/invariant 
and prove Theorems A and B 
(Theorem \ref{s2_Th1} and Corollary \ref{s2_Cor1}). 
In Section \ref{s3}, 
we define and study the modified scalar curvature. 
Then we review the necessary part of the $L^2$-index theory 
and prove Theorem~C (Theorem \ref{s3_T1}). 
Section \ref{s4} is devoted to the above examples 
and Theorem~D (Theorem \ref{s4_T1}). 

\noindent
{\bf \ref{int}.6. Acknowledgements.} 
Both authors are grateful 
to Harish Seshadri for interesting discussions. 
The first author is grateful to the Department of Mathematics 
at the University of Oregon for kind hospitality. 
The first author also would like to thank Tosiaki Kori and Mikio Furuta 
for useful comments on eta invariants and 
the Atiyah-Singer index theorem over orbifolds, respectively.

\section{Orbifold Yamabe constants/invariants}\label{s2}
{\bf \ref{s2}.1. Orbifolds with a finite number of singularities:
definitions.} There are several different approaches to and definitions
of orbifolds (cf.~\cite{Sa}, \cite{Kapovich}).  Since we work here with
orbifolds of a particular type, i.e., with only a finite number of
singularities, we give suitable definitions (cf.~\cite{Kr}). We
assume here that a finite group $\Gamma$ 
of the $n$-dimensional orthogonal group $O(n)$ 
inherits the standard $O(n)$-action on $\R^n$. 
We always assume that $n\geq 3$. 

First we define orbifolds with a finite number of singularities 
as an object in smooth category. 
\begin{Definition}\label{s2-d1}
{\rm Let $M$ be a locally compact Hausdorff space. We say that $M$ is
an $n$-dimensional \emph{orbifold with singularities} $$\Sigma_{\Gamma}
= \{(\check{p}_1,\Gamma_1), \ldots,(\check{p}_{s}, \Gamma_{s})\}$$ if the
following conditions are satisfied:
\begin{enumerate}
\item[{\bf (1)}] 
$\Sigma=\{\check{p}_1,\ldots,\check{p}_{s}\}\subset M$, 
and $M\setminus \Sigma$ is a smooth manifold of dimension $n$. 
\item[{\bf (2)}] 
$\Gamma = \{\Gamma_1, \ldots, \Gamma_s\}$ is the collection 
of subgroups $\Gamma_j$ of $O(n)$. 
Each group $\Gamma_j\ (j = 1, \ldots, s)$ is a nontrivial finite subgroup 
of $O(n)$ acting freely on $\R^n\setminus \{0\}$. 
\item[{\bf (3)}] 
For each $j\ (j = 1, \ldots, s)$, 
there exist an open neighborhood $U_j$ of $\check{p}_j$ 
and a homeomorphism $\phi_j: U_j \to \B_{\tau_j}(0)/\Gamma_j$ 
for some $\tau_j > 0 $ 
such that the restriction 
$
\phi_j : U_j \setminus \{\check{p}_j\} \to
\left(\B_{\tau_j}(0)/\Gamma_j\right)\setminus \{0\}
$ \ 
\ 
is a diffeomorphism. 
Here $\B_{\tau_j}(0) = 
\{x = (x^1, \ldots, x^n) \in \R^n\ |\ |x| < \tau_j \}$. 
\end{enumerate}
}
\end{Definition}
We refer to the pair $(\check{p}_j, \Gamma_j)$ 
as a \emph{singular point with the structure group} $\Gamma_j$ 
and the pair $(U_j,\phi_j)$ as a \emph{local uniformization}. 
Let $\pi_j : \B_{\tau_j}(0) \to \B_{\tau_j}(0)/\Gamma_j$ 
denote the canonical projection. 
\begin{Definition}\label{s2-d1-1}
{\rm Let $\widetilde{U}_j$ be an open neighborhood of $\check{p}_j$ 
and $\widetilde{\phi}_j : \widetilde{U}_j\to \B_{\tilde{\tau}_j}(0)/\Gamma_j$ 
a homeomorphism (for some $\tilde{\tau}_j>0$).  
We call the pair $(\widetilde{U}_j,\widetilde{\phi}_j)$ 
a \emph{compatible local uniformization} of $(\check{p}_j, \Gamma_j)$ 
if the following conditions are satisfied: 
\begin{enumerate}
\item[{\bf (1)}] 
The restriction 
$
\widetilde{\phi}_j : \widetilde{U}_j \setminus \{\check{p}_j\}
\to (\B_{\widetilde{\tau}_j}(0)\setminus \{0\})/\Gamma_j$ \ 
\ is a diffeomorphism. 
\item[{\bf (2)}] 
There exists a $\Gamma_j$-equivariant diffeomorphism 
$\Phi_j : \B_{\hat{\tau}_j}(0)\to \B_{\widetilde{\tau}_j}(0)$ 
for some $0<\hat{\tau}_j\leq \tau_j$ (possibly not onto) such that 
$$
\pi^*_j(\widetilde{\phi}_j\circ \phi_j^{-1})= \tilde{\pi}_j\circ \Phi_j \ \ 
\ \mbox{on} \ \ \B_{\hat{\tau}_j}(0)\setminus \{0\},
$$
where $\tilde{\pi}_j : \B_{\tilde{\tau_j}}(0) \to 
\B_{\tilde{\tau_j}}(0)/\Gamma_j$ denotes also the canonical projection.
\end{enumerate}}
\end{Definition}
To simplify the presentation, we assume, without particular mention, 
that an orbifold $M$ has only one singularity, i.e., 
$\Sigma_{\Gamma}=\{(\check{p},\Gamma)\}$. 
Let $\phi: U \to \B_{\tau}(0)/\Gamma$ be a local uniformization 
and $\pi: \B_{\tau}(0) \to \B_{\tau}(0)/\Gamma$ the canonical projection. 
We also always assume that $M$ is compact.  
Now we give the definition of orbifold metrics. 
\begin{Definition}\label{s2-d2}
{\rm A Riemannian metric $g\in \Riem(M\setminus \{\check{p}\})$ 
is an \emph{orbifold metric} 
if there exists a $\Gamma$-invariant smooth metric $\hat{g}$ on the disk $\B_{\tau}(0)$ 
such that $(\phi^{-1}\circ \pi)^*g = \hat{g}$ on $\B_{\tau}(0)\setminus \{0\}$.  }
\end{Definition}
We denote by $\Riem^{\orb}(M)$ the space of all orbifold metrics on
$M$.
\begin{Remark}
{\rm We note that if $g\in \Riem^{\orb}(M)$, 
then $R_{\hat{g}} \in C^{\infty}(\B_{\tau}(0))$, and hence 
the scalar curvature $R_g$ extends to 
a continuous function on $M$.} \hfill $\Box$ 
\end{Remark}
In the case when 
$\Sigma_{\Gamma} = 
\{(\check{p}_1,\Gamma_1),\ldots,(\check{p}_{s}, \Gamma_{s})\}$, 
the space $\Riem^{\orb}(M)$ of orbifold metrics is defined similarly. 

\noindent
{\bf \ref{s2}.2. Einstein orbifold metrics.} 
Let $M$ be a compact orbifold with a singularity 
$\Sigma_{\Gamma} = \{(\check{p},\Gamma)\}$. 
Consider the (normalized) Einstein-Hilbert functional 
$$ 
I : \Riem^{\orb}(M) \to \R, \ \ \ g\mapsto 
\frac{\int_{M} R_gd\sigma_g}{\Vol_g(M)^{\frac{n-2}{n}}}. 
$$ 
As in the case of smooth compact manifolds, 
we have the following fundamental result 
since the Stokes formula and the divergence theorem still hold over orbifolds. 
\begin{Proposition}\label{EH} 
The set of critical points of $I : \Riem^{\orb}(M) \to \R$ 
coincides with the set of Einstein orbifold metrics on $M$. 
\end{Proposition}

\noindent
{\bf \ref{s2}.3. Orbifold conformal classes.} 
We say that two orbifold metrics $g, 
\widetilde{g}\in \Riem^{\orb}(M)$ are \emph{pointwise conformal} 
if there exists a function 
$f \in C^0(M)\cap C^{\infty}(M\setminus \{\check{p}\})$ such that 
$$
\{
\begin{array}{l}
\widetilde{g} = e^{2f}\cdot g \ \ \mbox{on} \ \ M\setminus \{\check{p}\},
\\
(\phi^{-1}\circ \pi)^*f \in C^{\infty}(\B_\tau(0)).
\end{array}
\right.
$$
Here we used the composition of the maps 
$\B_{\tau}(0)\stackrel{\pi}\longrightarrow \B_\tau(0)/\Gamma 
\stackrel{\phi^{-1}}\longrightarrow U \stackrel{f}\longrightarrow \R$.  
Then we define an \emph{orbifold conformal class of $g$} as follows: 
$$ 
[g]_{\orb} := [g] \cap \Riem^{\orb}(M) 
= \{ e^{2f} \cdot g \ \left|
\begin{array}{l}
f\in C^0(M)\cap C^{\infty}(M\setminus\{\check{p}\}),
\\
(\phi^{-1}\circ \pi)^*f  \in C^{\infty}(\B_\tau(0))
\end{array}
\right.\} . 
$$ 
Let ${\mathcal C}^{\orb}(M)$ 
denote the space of all orbifold conformal classes.  
The proofs of the following two lemmas are similar 
to the case of smooth compact manifolds (see \cite{Aubin,KW}). 
\begin{Lemma}\label{lemma1-5}
Let $g\in\Riem^{\orb}(M)$ be a metric satisfying 
one of the following conditions: 
$R_g>0$, $R_g\equiv 0$ or $R_g< 0$ everywhere on $M$. 
Then for another such metric $\tilde{g}\in[g]_{\orb}$, 
the sign of $R_{\tilde{g}}$ is identical with the sign of $R_g$. 
\end{Lemma}
\begin{Lemma}\label{lemma1-6}
Let $g\in\Riem^{\orb}(M)$ be a metric with $R_g\equiv\const\leq 0$. 
Then for any metric $\tilde{g}\in [g]_{\orb}$ with 
$R_{\tilde{g}}\equiv\const$, there exists a constant $k>0$ such that
$\tilde{g}=k\cdot g$.
\end{Lemma}
An orbifold $M$ is called \emph{good} if its universal cover 
$\widetilde{M}$ is a smooth manifold (cf.~\cite{Kapovich}). 
Let $\pi: \widetilde{M}\to M$ be the projection. 
\begin{Proposition}
Let $g\in\Riem^{\orb}(M)$ be an Einstein metric with
$R_g\equiv\const>0$. Then any metric $\tilde{g}\in [g]_{\orb}$ with
$R_{\tilde{g}}\equiv\const$ is also Einstein. Furthermore, if $M$ is a
good orbifold, then $\tilde{g}=k\cdot g$ for a constant $k>0$ unless
$(\widetilde{M},\pi^*g)$ is isometric to $(S^n,g_S)$, 
where $g_S \in \Riem(S^n)$ denotes the standard metric 
of constant curvature one. 
\end{Proposition}
\begin{Proof}
The first statement follows directly from 
the proof of Proposition 1.4 in \cite{Schoen} 
combined with the divergence theorem for orbifolds 
(see also \cite[Section 2]{Escobar}). 

Now assume that the universal cover $\widetilde{M}$ is a smooth manifold. 
Let $\pi: \widetilde{M}\to M$ denote the projection. 
Let $\widetilde{g} = u^{\frac{4}{n-2}}\cdot g$ 
be an orbifold metric with $R_{\widetilde{g}} \equiv \textrm{const.}$ on $M$. 
Then unless $(\widetilde{M},\pi^*g)\cong (S^n,g_S)$, 
\cite[Theorem~A]{Obata} implies that 
$\pi^*u\equiv \const>0$ on $\widetilde{M}$, 
and hence $u\equiv \const>0$ on $M$. 
\end{Proof}

\noindent
{\bf \ref{s2}.4. Orbifold Yamabe constants/invariants.} 
We define the \emph{orbifold Yamabe constant} $Y_{[g]_{\orb}}(M)$ 
of $(M,[g]_{\orb})$ as follows:
$$ 
Y_{[g]_{\orb}}(M) := \inf_{\widetilde{g}\in [g]_{\orb}} 
\frac{\int_M R_{\widetilde{g}} d\sigma_{\widetilde{g}}} 
{\Vol_{\widetilde{g}}(M)^{\frac{n-2}{n}}}\ . 
$$
Similarly to the case of smooth compact manifolds, 
we have Aubin's inequality 
$-\infty < Y_{[g]_{\orb}}(M) \leq Y(S^n)$ (cf.~\cite{Aubin}). 
Then we define the \emph{orbifold Yamabe invariant of $M$} 
$$ 
Y^{\orb}(M) := \sup_{C\in {\mathcal C}^{\orb}(M)} Y_{C}(M) \ \ \ (\ 
\leq Y(S^n)\ ). 
$$ 
Now we state some technical facts without proofs. 
The following result also follows from 
the divergence theorem for orbifolds. 
\begin{Lemma}\label{s2_L1} 
For an orbifold metric $\widetilde{g}\in [g]_{\orb}$ 
with $\widetilde{g}=u^{\frac{4}{n-2}}g $, then 
$$
\int_M R_{\widetilde{g}} d\sigma_{\widetilde{g}} =
\int_M\left(\alpha_n|du|^2_g +R_gu^2\right) d\sigma_g ,
\ \ \ \ \mbox{where \ \  $\alpha_n=\frac{4(n-1)}{n-2}>0$}.
$$
\end{Lemma}

We consider the Yamabe functional  
$$
Q_g(u):= \frac{E_g(u)}{\left(\int_M |u|^{\frac{2n}{n-2}} d\sigma_g 
\right)^{\frac{n-2}{n}}}=
\frac{\int_M\left(\alpha_n|du|^2_g +R_gu^2 \right) 
d\sigma_g}{\left(\int_M |u|^{\frac{2n}{n-2}} d\sigma_g 
\right)^{\frac{n-2}{n}}} 
$$
for any $u\in C^0(M)\cap C^{\infty}(M\setminus \{p\})$ 
with $u\not\equiv 0$ and $(\phi^{-1}\circ \pi)^*u\in C^{\infty}(\B_\tau(0))$. 

\noindent
\begin{Remark}{\rm We recall that for a compact smooth manifold $N$ 
of dim $N \geq 3$, any point $q\in N$ and any metric $h \in \Riem(N)$, 
the Yamabe constant $Y_{[h]}(N)$ satisfies} 
$$ 
Y_{[h]}(N)= \inf_{\begin{array}{c}^{u\in  C_c^{\infty}(N \setminus \{q\})}\\ 
^{u\not\equiv 0}\end{array}} Q_h (u). 
$$ 
\end{Remark}

\noindent
The orbifold Yamabe constant has a similar property. 
\begin{Lemma}\label{s2_L2} 
The orbifold Yamabe constant $Y_{[g]_{\orb}}(M)$ satisfies 
$$ 
Y_{[g]_{\orb}}(M) = 
\inf_{\begin{array}{c}^{u\in  C_c^{\infty}(M\setminus \{\check{p}\})}\\ 
^{u\not\equiv 0}\end{array}} Q_g(u)\ . 
$$ 
Here $\check{p}$ is the singular point of $M$. 
\end{Lemma}
{\bf \ref{s2}.5. Approximation for orbifold metrics.} 
At the singular point $\check{p} \in U (\subset M)$, 
we prove the following approximation, 
which is an orbifold version of Kobayashi's approximation lemma 
(\cite[Lemma~3.2]{Ko}). 
This approximation is a crucial tool 
linking the orbifold Yamabe invariant to 
the corresponding $h$-cylindrical Yamabe invariant. 
\begin{Proposition}\label{s2_P1}
Let $g \in \Riem^{\orb}(M)$ be an orbifold metric and 
$\hat{g} \in \Riem(\bold{B}_{\tau}(0))$ the lifting metric with 
$\hat{g} = (\varphi^{-1} \circ \pi)^*g$ 
on $\bold{B}_{\tau}(0) \setminus \{0\}.$ 
Let $y = (y^1, \ldots, y^n)$ denote normal coordinates 
with respect to $\hat{g}$ on an open neighborhood 
$\bold{U}\ (\subset \bold{B}_{\tau}(0))$ of $0$, 
and $\bold{U}_{\rho} = 
\{y \in \bold{U}\ |\ \text{\rm dist}_{\hat{g}}(0, y) = |y| < \rho \}$ 
for small $\rho > 0$. 
Then for any small 
$\delta>0$ there exists a metric $g_{\delta}\in \Riem^{\orb}(M)$ 
such that: 
\begin{enumerate}
\item[{\bf (1)}] 
$g_{\delta} = g$ on $M\setminus V_{\delta}$, 
where 
$V_{\delta} = \{ q \in M\ |\ \text{\rm dist}_g(\check{p}, q) < \delta \}$, 
\item[{\bf (2)}] 
$\hat{g}_{\delta}= (\phi^{-1}\circ\pi)^*g_{\delta}$ 
is pointwise conformal to the flat metric $\sum_{i=1}^n(dy^i)^2$ 
on $\bold{U}_{\epsilon(\delta)}(0) \setminus \{0\}$ 
for a constant $\epsilon(\delta)\ (0 < \epsilon(\delta) < \delta)$, 
where $\hat{g}_{\delta} \in \Riem(\bold{B}_{\tau}(0))$ 
denotes the corresponding lifting metric to $g_{\delta}$, 
\item[{\bf (3)}] 
$g_{\delta} \to g$ uniformly on $M\setminus \{\check{p}\}$ 
and $\hat{g}_{\delta}\to \hat{g}$ uniformly on $\B_{\tau}(0)$ 
as $\delta\to 0$, 
\item[{\bf (4)}] 
$R_{g_{\delta}} \to R_g$ uniformly on $M$ as $\delta\to 0$. 
\end{enumerate}
\end{Proposition}
\begin{Proof}
Let $g_0=\sum_{i=1}^n(dy^i)^2$ be the flat metric on $\bold{U}$. 
Then 
$$
\hat{g}(y) = g_0(y) + O(|y|^2) \ \ \ \mbox{on } \ \ \bold{U}. 
$$
Since $\hat{g}$ is $\Gamma$-invariant, the exponential map 
(with respect to $\hat{g}$) commutes with $\Gamma$, 
and hence the flat metric $g_0$ and the distance function 
$\text{\rm dist}_{\hat{g}}(0, y) = |y|$ from $0$ to $y$ 
are also $\Gamma$-invariant. 
Moreover, $(\bold{U}, y = (y^1, \ldots, y^n))$ 
is a compatible local uniformization. 
We denote $\kappa:=R_{\hat{g}}(0)$, and consider the metric 
$$ 
\tilde{g} := 
\left(1-\frac{\kappa}{2\alpha_n}|y|^2\right)^{\frac{4}{n-2}}\cdot 
g_0\ \ \ \ \mbox{on } \ \ \ \bold{U}. 
$$ 
Then it follows that $R_{\tilde{g}}(0)=\kappa= R_{\hat{g}}(0)$. Now we
use the cut-off function $w_{\delta}$ (given in \cite{Ko}) to 
construct the following approximation $\hat{g}_{\delta}$:
$$
\hat{g}_{\delta}:=\hat{g}+ w_{\delta}(r)(\tilde{g}-\hat{g}) 
\ \ \ \mbox{on } \ \ \B_\tau(0), \ \ \ r=|y|.
$$
Since $j_0^1(\tilde{g})=j_0^1(\hat{g})$ and 
$R_{\tilde{g}}(0)=R_{\hat{g}}(0)$, Kobayashi's approximation technique 
\cite{Ko} implies that $\hat{g}_{\delta}$ satisfies the above conditions 
{\bf (1)}--{\bf (4)} on $\B_\tau(0)$. 
By the construction, the metrics $\hat{g}$, 
$\tilde{g}$ and the function $w_{\delta}$ are $\Gamma$-invariant, 
and hence the metric $\hat{g}_{\delta}$ is also $\Gamma$-invariant. 
This implies that there exists an orbifold metric 
$g_{\delta}\in\Riem^{\orb}(M)$ such that 
$\hat{g}_{\delta} = (\phi^{-1}\circ\pi)^*g_{\delta}$. 
Then the metric $g_{\delta}$ 
satisfies the above conditions {\bf (1)}--{\bf (4)}.
\end{Proof} 
The following lemma is an analogue of the case of smooth compact manifolds 
(cf.~\cite{Besse, Ko}). 
\begin{Lemma}\label{s2_L4} 
Let $g_{\delta}, \ g\in \Riem^{\orb}(M)$ be metrics satisfying 
$$
g_{\delta}\to g \ \ \mbox{and}\ \ \ R_{g_{\delta}}\to R_g \ \ \
\mbox{uniformly on $M$ as $\delta\to 0$}. 
$$
Then $Y_{[g_{\delta}]_{\orb}}(M)\to Y_{[g]_{\orb}}(M)$ as $\delta\to 0$.
\end{Lemma}
{\bf \ref{s2}.6. Orbifold and $h$-cylindrical Yamabe invariants.} 
We recall briefly the following: 
The open manifold $M\setminus \{\check{p}\}$ 
equipped with an appropriate cylindrical metric can be considered 
as a cylindrical manifold (see \cite{AB2} or Section \ref{s3}). 
For a complete metric $\bar{g}\in \Riem(M \setminus \{\check{p}\})$, 
the open manifold $(M\setminus \{\check{p}\},\bar{g})$ is called 
a \emph{cylindrical manifold modeled by} $(S^{n-1}/\Gamma, h_{\Gamma})$ 
if there exist an open neighborhood $V (\subset U)$ of $\check{p}$ and 
a coordinate system $(z,t) \in (S^{n-1}/\Gamma) \times [0,\infty) 
\cong V \setminus \{\check{p}\}$ such that 
$$
\bar{g}(z,t)=h_{\Gamma}(z)+ dt^2 \ \ \ \ \mbox{on} \ \ 
(S^{n-1}/\Gamma) \times [1,\infty).
$$
The \emph{cylindrical Yamabe constant}
$Y^{\cyl}_{[\bar{g}]}(M\setminus \{\check{p}\})$ is defined by
$$
Y^{\cyl}_{[\bar{g}]}(M\setminus \{\check{p}\}) := 
\inf_{\begin{array}{c}^{u\in  C_c^{\infty}(M\setminus \{\check{p}\})}\\
^{u\not\equiv 0}\end{array}} Q_{\bar{g}}(u) \ \ \ (\ \leq  Y(S^n)\ ).
$$
\vspace{-7mm}

Let $\Riem^{h_{\Gamma}\hbox{-}\cyl}(M \setminus \{\check{p}\})\ 
(\subset \Riem(M\setminus \{\check{p}\}))$ 
be the space of cylindrical metrics modeled by 
$(S^{n-1}/\Gamma, h_{\Gamma})$. 
Then the \emph{$h_{\Gamma}$-cylindrical Yamabe invariant} 
$Y^{h_{\Gamma}\hbox{-}\cyl}(M \setminus \{\check{p}\})$ is also defined by 
$$
Y^{h_{\Gamma}\hbox{-}\cyl}(M\setminus \{\check{p}\}):= 
\sup_{\bar{g} \in \Riem^{h_{\Gamma}\hbox{-}\cyl}(M\setminus \{\check{p}\})} 
Y^{\cyl}_{[\bar{g}]}(M\setminus \{\check{p}\}) \ \ \ (\ \leq  Y(S^n)\ ). 
$$
\begin{Theorem}\label{s2_Th1}
Let $M$ be a compact orbifold with a singularity 
$\Sigma_{\Gamma} = \{(\check{p}, \Gamma)\}$. 
Then $Y^{\orb}(M)=Y^{h_{\Gamma}\hbox{-}\cyl}(M \setminus \{\check{p}\})$. 
\end{Theorem}
Theorem \ref{s2_Th1} and \cite[Propositions 2.11, 2.12, 6.5]{AB2} imply the
following result.
\begin{Corollary}\label{s2_Cor1}\qquad 
Let $M$ be a compact orbifold with singularities 
$\Sigma_{\Gamma} =$ 
$\{(\check{p}_1,\Gamma_1), \ldots,(\check{p}_{s}, \Gamma_{s})\}$. 
Then 
$$
Y^{\orb}(M) = 
Y^{h_{\Gamma}\hbox{-}\cyl}(M\setminus \{\check{p}_1,\ldots,\check{p}_{s}\}) 
\leq \min_{1\leq j\leq s}\frac{Y(S^n)}{|\Gamma_j|^{\frac{2}{n}}}.
$$
\end{Corollary}
\begin{lproof}{Proof of Theorem \ref{s2_Th1}} 
First we prove that 
$Y^{\orb}(M) \leq Y^{h_{\Gamma}\hbox{-}\cyl}(M\setminus \{\check{p}\})$. 
By Proposition \ref{s2_P1} and Lemma \ref{s2_L4}, 
for any $\epsilon>0$ and any metric $g \in \Riem^{\orb}(M)$ 
there exists a metric $g_{\epsilon} \in \Riem^{\orb}(M)$ such that 
$$
\{\!\!
\begin{array}{l}
|Y_{[g]_{\orb}}(M)-Y_{[g_{\epsilon}]_{\orb}}(M)|<\epsilon, \\
\displaystyle
(\phi^{-1}\!\circ\!\pi)^*g_{\epsilon}\ \simeq\ g_0\!=\!
\sum_{i=1}^n(dx^i)^2 \ 
\mbox{: pointwise conformal near} \ 0\in \B_\tau(0).
\end{array}
\right.
$$
Let $\bar{g}_{\epsilon}\in [g]_{\orb}$ be a cylindrical metric 
on $M \setminus \{\check{p}\}$ 
with $\bar{g}_{\epsilon}(z,t)=h_{\Gamma}(z)+dt^2$ 
on the end 
$(S^{n-1}/\Gamma)\times [1, \infty) \cong V \setminus \{\check{p}\}$. 
Then we use Lemma \ref{s2_L2} to show 
$$
\begin{array}{rcl}
Y_{[g_{\epsilon}]_{\orb}}(M) &= &
\displaystyle
\inf_{\begin{array}{c}^{u\in  C_c^{\infty}(M\setminus \{\check{p}\})}\\
^{u\not\equiv 0}\end{array}} Q_{g_{\epsilon}}(u)
= 
\displaystyle 
\inf_{\begin{array}{c}^{u\in  C_c^{\infty}(M\setminus \{\check{p}\})}\\
^{u\not\equiv 0}\end{array}} Q_{\bar{g}_{\epsilon}}(u)
\\
\\
 &= &
Y^{\cyl}_{[\bar{g}_{\epsilon}]}(M \setminus \{\check{p}\}) \leq 
Y^{h_{\Gamma}\hbox{-}\cyl}(M \setminus \{\check{p}\}).
\end{array}
$$
This implies that 
$Y^{\orb}(M) \leq Y^{h_{\Gamma}\hbox{-}\cyl}(M\setminus\{\check{p}\})$.  

Second we prove that 
$Y^{\orb}(M) \geq Y^{h_{\Gamma} \hbox{-}\cyl}(M\setminus\{\check{p}\})$. 
We start with an arbitrary cylindrical metric 
$\bar{g} \in \Riem^{h_{\Gamma}\hbox{-}\cyl}(M\setminus \{\check{p}\})$ 
and a cylindrical coordinate system 
$(z,t) \in (S^{n-1}/\Gamma) \times [1,\infty) \cong 
V \setminus \{\check{p}\}$ as above. 

Set $r=e^{-t}$, then $ \bar{g}(z,t)=r^{-2}(dr^2 +r^2\cdot h_{\Gamma}(z))$ 
on  the cylinder $(S^{n-1}/\Gamma) \times [1,\infty)$. 
Note that on 
$(S^{n-1}/\Gamma) \times [1,\infty)$ the metric 
$r^2\cdot \bar{g}(z,t)= dr^2 +r^2\cdot h_{\Gamma}(z)$ 
is flat and is extended to a smooth metric 
$g_0$ on $M\setminus \{\check{p}\}$. 
This implies that there exists a homeomorphism 
$\tilde{\phi}: \tilde{U} \to 
\B_{\tilde{\tau}}(0)/\Gamma$ for a constant $\tilde{\tau}>0$ such that 
\begin{enumerate}
\item[{\bf (1)}] $\tilde{\phi}: (\tilde{U} \setminus \{\check{p}\})\to
(\B_{\tilde{\tau}}(0)\setminus \{0\})/\Gamma$ is a diffeomorphism,
\item[{\bf (2)}] 
$(\tilde{\phi}^{-1}\circ \tilde{\pi})^*g_0 = \sum_{i=1}^n(dy^i)^2$, 
where $y = (y^1, \ldots, y^n) \in \tilde{U}$ 
are normal coordinates around $\check{p}$ with respect to $g_0$, 
\item[{\bf (3)}] the diffeomorphism $\tilde{\phi}\circ\phi^{-1}$ is
lifted to a $\Gamma$-equivariant diffeomorphism $\Phi:
\B_{\hat{\tau}}(0)\to \B_{\tilde{\tau}}(0)$ onto its image 
for a constant $0< \hat{\tau} \leq {\tau}$.
\end{enumerate}
In other words, $(\tilde{U}, \tilde{\phi})$ 
is a compatible local uniformization, 
and the above metric $g_0$ 
is an orbifold metric on $M$. 
>From Lemma \ref{s2_L2}, we obtain that 
$$
\begin{array}{l}
Y_{[g_0]_{\orb}}(M) = 
\displaystyle
\!\!\!\!\!\!
\inf_{\begin{array}{c}^{u\in  C_c^{\infty}(M\setminus \{\check{p}\})}\\
^{u\not\equiv 0}\end{array}}\!\!\! Q_{g_0}(u)
= 
\displaystyle
\!\!\!\!\!\!
\inf_{\begin{array}{c}^{u\in  C_c^{\infty}(M\setminus \{\check{p}\})}\\
^{u\not\equiv 0}\end{array}}\!\!\! Q_{\bar{g}}(u)
=
Y^{\cyl}_{[\bar{g}]}(M\setminus \{\check{p}\}) .
\end{array}
$$
This implies that 
$Y^{\orb}(M) \geq Y^{h_{\Gamma}\hbox{-}\cyl}(M\setminus \{\check{p}\})$. 
\end{lproof}
\section{$L^2$-harmonic spinors and Yamabe invariant}\label{s3}
{\bf \ref{s3}.1. Modified scalar curvature.} 
Throughout the rest of the paper, we consider only $4$-manifolds.  
Let $X$ be an open smooth $4$-manifold with tame ends, i.e., 
it is diffeomorphic to $W \cup_Z (Z \times [0,\infty))$, 
where $W (\subset X)$ is a compact submanifold of dim $W = 4$ 
with boundary $\p W \cong Z\times \{0\}$ (possibly disconnected), 
see Fig. \ref{s3}.1.

\hspace*{2mm}\PSbox{orb2a.pstex}{35mm}{35mm}
\begin{picture}(5,0)
\put(-5,20){{\small $W$}}
\put(40,5){{\small $Z\times [0,\infty)$}}
\put(10,45){{\small $Z\times \{0\}$}}
\put(35,40){{\small $Z\times \{1\}$}}
\end{picture}
\vspace{3mm}

\centerline{{\small {\bf Fig. \ref{s3}.1.} Cylindrical manifold
$(X,\bar{g})$.}}
\vspace{2mm}

\noindent 
Let $\Riem^{\cyl}(X)$ denote the space of all cylindrical metrics on $X$. 
We choose a cylindrical metric $\bar{g} \in \Riem^{\cyl}(X)$. 
In terms of \cite{AB2}, $(X, \bar{g})$ is a 
\emph{cylindrical $4$-manifold modeled by} $(Z,h)$ for $h \in \Riem(Z)$ 
if there exists a product coordinate system $(z,t)$ 
on the end $Z \times [0, \infty)$ such that 
$\bar{g}(z,t) = h(z) + dt^2$ on $Z \times [1, \infty)$. 
Below we use the notation that $\p_{\infty}\bar{g} = h$.  

Let $\omega\in \Omega^2(X)$ be a 2-form. 
We denote by $\L_{\bar{g}}=-6\Delta_{\bar{g}}+ R_{\bar{g}}$ 
the conformal Laplacian of $\bar{g}$.  
We define the \emph{modified conformal Laplacian 
$\L_{(\bar{g},\omega)}$ of $(\bar{g},\omega)$} by
$$
\L_{(\bar{g},\omega)}:=\L_{\bar{g}}-|\omega|_{\bar{g}} =
 -6\Delta_{\bar{g}}+ (R_{\bar{g}}-|\omega|_{\bar{g}}).
$$
Note that the pointwise norm $|\omega|_{\bar{g}}$ 
is only of $C^{1,0}$ class on $X$, in general. 

\begin{Fact}\label{s3_F1}{\rm (cf.~\cite{GL})} 
Set $\check{g}=u^2\cdot\bar{g}$ for $u\in C^{\infty}_+(X)$.
Then $\L_{\check{g}}(f) = u^{-3} \cdot \L_{\bar{g}}(uf)$ 
for $f\in C^2(X)$. 
\end{Fact}
\begin{Definition}{\rm (Gursky-LeBrun \cite{Gu,GL}) 
For $\bar{g}\in\Riem^{\cyl}(X)$, 
let $\check{g}=u^2\cdot\bar{g}$ be a conformal metric, 
where $u\in C^{\infty}_+(X)$. 
We call the function 
$R_{(\check{g},\omega)}:=R_{\check{g}} -|\omega|_{\check{g}}$ 
the \emph{modified scalar curvature of $(\check{g},\omega)$}.}
\end{Definition}
Set $E_{(\bar{g},\omega)}(f):=\int_X\left[
6|df|^2_{\bar{g}} + (R_{\bar{g}}-|\omega|_{\bar{g}})f^2 
\right]d\sigma_{\bar{g}}$ for $f\in L^{1,2}_{\bar{g}}(X)$, where
$L^{1,2}_{\bar{g}}(X)$ denotes the Sobolev space of square-integrable
functions on $X$ (with respect to $\bar{g}$) up to their first weak
derivatives.  Then we consider the functional
$$
\mathbb{Q}_{(\bar{g},\omega)}(f):= 
\frac{E_{(\bar{g},\omega)}(f)} {\int_X f^2d\sigma_{\bar{g}}} 
\qquad \textrm{for}\quad f \in L^{1,2}_{\bar{g}}(X)\ 
\textrm{with}\ f\not\equiv 0 
$$
and the bottom of the spectrum of the operator
$\L_{(\bar{g},\omega)}$
$$
\lambda_{(\bar{g},\omega)}:= 
\inf_{\begin{array}{c}^{f\in L^{1,2}_{\bar{g}}(X)}\\
^{f\not\equiv 0}\end{array}} \mathbb{Q}_{(\bar{g},\omega)}(f).
$$
>From now on we assume that the $2$-form $\omega$ has a compact support. 
\begin{Lemma}\label{s3_L1}
For $\bar{g}\in \Riem^{\cyl}(X)$ and $\omega\in \Omega_{\cp}^2(X)$ as above, 
$\lambda_{(\bar{g},\omega)}>-\infty$. 
\end{Lemma}
\begin{Proof}
Indeed, $\lambda_{(\bar{g},\omega)}\geq \displaystyle\inf_X
(R_{\bar{g}}-|\omega|_{\bar{g}})> -\infty$ since $\bar{g}\in
\Riem^{\cyl}(X)$ and $\omega\in \Omega_{\cp}^2(X)$.
\end{Proof}
On the slice manifold $(Z,h)$, we consider the \emph{almost conformal
Laplacian} ${\mathcal L}_h:=-6\Delta_h+ R_h$ (see \cite{AB2}). 
\vspace{2mm}

\noindent
{\bf Convention 1.}  From now on we always assume that 
$\lambda_h:=\lambda({\mathcal L}_h)>0$, 
where $\lambda({\mathcal L}_h)$ denotes the first eigenvalue 
of ${\mathcal L}_h$. 
\begin{Remark}{\rm (\cite[Section 2]{AB2})
It is easy to see that the condition $R_h>0$ everywhere on $Z$
provides a sufficient condition for $\lambda_h>0$. 
On the other hand,
if $\lambda_h\leq 0$, then $Y^{\cyl}_{[\bar{g}]}(X)\leq 0$. 
Moreover, if $\lambda_h<0$, 
then $Y^{\cyl}_{[\bar{g}]}(X)=-\infty$.}
\end{Remark}
\begin{Lemma}\label{s3_L2}
Assume that $\lambda_{(\bar{g},\omega)}<\lambda_h$. 
Then there exists a function 
$u\in C^{2,\alpha}_+(X)\cap L_{\bar{g}}^{1,2}(X)$ (for any $0<\alpha<1$) 
such that $\L_{(\bar{g},\omega)}u = \lambda_{(\bar{g},\omega)}u$. 
\end{Lemma}
{\em Proof.}
Consider the cylinder $Z\times [\ell,\infty)$ with $\ell\geq 1$ and 
denote $X(\ell):= W\cup_{Z} (Z\times [0,\ell])$.  
Note that if there exists a constant $\ell\geq 1$ satisfying
$$
\lambda(\ell) := 
\inf_{\begin{array}{c}^{f\in C_c^{\infty}(Z\times [\ell,\infty))}\\
^{f\not\equiv 0}\end{array}} \mathbb{Q}_{(\bar{g},\omega)}(f) > 
\lambda_{(\bar{g},\omega)},
$$
then the standard argument implies the existence of a non-zero minimizer 
$u\in C^{2,\alpha}_+(X)\cap L_{\bar{g}}^{1,2}(X)$ 
of the functional $\mathbb{Q}_{(\bar{g},\omega)}$. 
Hence it is enough to prove that 
$\lambda(\ell) \geq \lambda_h$ 
for a constant $\ell\geq 1$. 
We choose $\ell$ satisfying $\supp(\omega)\subset X(\ell)$. 
Then $R_{\bar{g}}=R_h$ and 
$|\omega|_{\bar{g}}=0$ on $Z\times [\ell,\infty)$. Hence for any $f\in
C^{\infty}_c(Z\times [\ell,\infty))$,
$$
\begin{array}{rcl}
\displaystyle
\left.E_{(\bar{g},\omega)}\right|_{Z\times [\ell,\infty)}(f) \!\!&=&\!\!
\displaystyle
\int_{Z\times [\ell,\infty)}\left( 6|df|^2 + (R_{\bar{g}} -
|\omega|_{\bar{g}})f^2\right) d\sigma_{\bar{g}} 
\\
\\
\!\!&=&\!\!
\displaystyle
\int_{[\ell,\infty)} dt \int_{Z}\left[
6(\p_t f)^2 + 6 |\nabla^Z f|^2_h + R_h f^2
\right]d\sigma_h 
\\
\\
\!\!&\geq&\!\!
\displaystyle
\int_{[\ell,\infty)}\!\! 
dt \left(\! \lambda_h \cdot \int_{Z} f^2 d\sigma_h\! \right)
\geq \lambda_h \cdot 
\int_{Z\times [\ell,\infty)}f^2 d\sigma_{\bar{g}} . 
\ \Box
\end{array}
$$
The next lemma shows that we may always assume the condition 
of Lemma \ref{s3_L2} without loss of generality. 
The technique below of changing the given metric conformally 
within a compact set is known as the 
\emph{conformal-rescaling argument} 
(see also the proofs of \cite[Theorem~3]{AB1} and 
\cite[Proposition~7.3]{AB2}). 
\begin{Lemma}\label{s3_L3} 
There exists a metric $\tilde{g}\in [\bar{g}]$ 
with $\tilde{g}\equiv \bar{g}$ on $X\setminus X(1)$ such that 
$\lambda_{(\tilde{g},\omega)}<\lambda_h$. 
In particular, $\tilde{g} \in \Riem^{\cyl}(X)$. 
\end{Lemma}
\begin{Proof} 
When $\lambda_{(\bar{g}, \omega)} \leq 0$, 
the inequality $\lambda_{(\bar{g}, \omega)} < \lambda_h$ is trivial. 
Hence we consider only the case when 
$\lambda_{(\bar{g}, \omega)} > 0$. 
We change the cylindrical metric $\bar{g}$ 
to another cylindrical metric 
$\bar{g}_v = e^{2v}\cdot \bar{g}$, 
where $v \in C^{\infty}(X)$ with 
$$
v = \{
\begin{array}{l}
k\equiv \const>0 \ \ \mbox{on} \ \ W= X\setminus (Z\times (0,\infty)),
\\
0 \ \  \mbox{on} \ \ Z\times [1,\infty).
\end{array}
\right.
$$
In particular, 
$\L_{(\bar{g}_v, \omega)} = e^{-2k}\cdot \L_{(\bar{g}_v, \omega)}$ on $W$. 
We choose $k >> 1$ sufficiently large. 
Then the Dirichlet first eigenvalue 
$\lambda_{(\overline{g}_v, \omega)}(W)$ on $W$ satisfies that 
$\lambda_{(\overline{g}_v, \omega)}(W) < \lambda_h$.  
Hence by the domain
monotonicity of the Dirichlet eigenvalues, 
we then obtain that 
$\lambda_{(\bar{g}_v, \omega)} \leq 
\lambda_{(\bar{g}_v, \omega)}(W) < \lambda_h$.  
To complete the proof, 
we let $\tilde{g} = \bar{g}_v$.
\end{Proof}

\noindent
{\bf Convention 2.} From now on we always assume that 
$\lambda_{(\bar{g}, \omega)} < \lambda_h$ since 
$Y^{\cyl}_{[\tilde{g}]}(X)$ $=Y^{\cyl}_{[\bar{g}]}(X)$.
\vspace{2mm}

\noindent
Lemmas \ref{s3_L2}, \ref{s3_L3} imply the following lemma, 
whose proof is similar 
to the case of compact smooth manifolds (see \cite{GL}). 
\begin{Lemma}\label{s3_L4}
There exists a function 
$u\in C_+^{2,\alpha}(X)\cap L_{\bar{g}}^{1,2}(X)$ (where $0<\alpha<1$) 
such that for the metric $\check{g}=u^2\cdot \bar{g}$ 
$$
R_{(\check{g},\omega)} \ \ \{
\begin{array}{l}
>0 
\\
\equiv 0
\\
<0 
\end{array}
\right. \ \ \ \mbox{everywhere on $X$}.
$$
Furthermore, these cases are mutually exclusive.
\end{Lemma}
The next result follows from \cite[Proposition 7.1]{AB2}.
\begin{Fact}\label{s3_F2}
Let $\bar{g},\tilde{g}\in \bar{C}\cap\Riem^{\cyl}(X)$ 
be two cylindrical metrics 
which are pointwise conformal. Then 
$\sign(\lambda_{(\bar{g},\omega)})=
\sign(\lambda_{(\tilde{g},\omega)})$.
\end{Fact}
Now for $\bar{g}\in \Riem^{\cyl}(X)$, 
let $[\bar{g}]_{L^{1,2}_{\bar{g}}}\ (\subset [\bar{g}])$ denote 
the $L^{1,2}_{\bar{g}}$-\emph{conformal class} of $\bar{g}$ 
consisting of all metrics $u^2 \cdot \bar{g}$, 
where $u \in C^{\infty}_+(X) \cap L^{1,2}_{\bar{g}}(X)$ (see \cite{AB2}). 
Note that the $L^2_{\bar{g}}$-norm 
$
\|\zeta\|_{L^2_{\bar{g}}}:=\left(
\int_X |\zeta|^2_{\bar{g}} d\sigma_{\bar{g}}\right)^{\frac{1}{2}}$ 
of $\zeta\in \Omega^2(X)$ 
depends only on the conformal class $[\bar{g}]$, 
and hence set 
$\|\zeta\|_{L^2_{[\bar{g}]}}:=\|\zeta\|_{L^2_{\bar{g}}}$.  
Then, Fact \ref{s3_F2} combined with the argument 
given in \cite[Corollary~4]{GL} implies the following assertion. 
\begin{Proposition}\label{s3_P1}
Let $\bar{g}$ be a cylindrical metric on $X$, 
and $\omega\in \Omega^2_{\cp}(X)$ as above. 
Then one of the following {\bf (1)--(3)} holds: 
\begin{enumerate}
\item[{\bf (1)}] 
There exists a metric $\check{g} \in [\bar{g}]_{L^{1,2}_{\bar{g}}}$ 
such that $R_{\check{g}}> |\omega|_{\check{g}}$ on $X$, 
\item[{\bf (2)}] 
$Y^{\cyl}_{\bar{C}}(X) < \|\omega\|_{L^2_{\bar{C}}}^2$, 
\item[{\bf (3)}] 
$Y^{\cyl}_{\bar{C}}(X) = \|\omega\|_{L^2_{\bar{C}}}= 0$, 
and there exists a metric 
$\check{g}\in [\bar{g}]_{L^{1,2}_{\bar{g}}}$ such that 
$R_{\check{g}} \equiv 0$ on $X$. 
\end{enumerate} 
\end{Proposition}
{\bf \ref{s3}.2. $L^2$-harmonic 2-forms.} 
Let $(X,\bar{g})$ be a cylindrical $4$-manifold modeled by $(Z,h)$. 
A $2$-form $\zeta \in \Omega^2(X)$ is 
an $L^2_{\bar{g}}$-harmonic $2$-form on $X$, i.e., $\zeta$ satisfies 
$$
\{
\begin{array}{l}
d\zeta = d^{*}\zeta = 0,
\\
\|\zeta\|^2_{L^2_{\bar{g}}}:= 
\int_X|\zeta|^2_{\bar{g}} d\sigma_{\bar{g}} <\infty.
\end{array}
\right.
$$
Here $ d^{*}$ stands for the codifferential 
with respect to the metric $\bar{g}$. 
We consider the space of $L^2_{\bar{g}}$-harmonic $2$-forms 
$$
{\mathcal H}_{\bar{g}}^2(X):= \{ \zeta\in \Omega^2(X) \ 
| \ \zeta \ \ \mbox{is $L^2_{\bar{g}}$-harmonic}\ \} .
$$
We recall the following well-known facts.
\begin{Fact}\label{s3_F3}
{\rm (cf.~\cite{Do2})} 
Let $\zeta\in {\mathcal H}_{\bar{g}}^2(X)$ 
be an $L^2_{\bar{g}}$-harmonic 2-form. 
Then the function $|\zeta|_{\bar{g}}$ 
decays exponentially on $Z\times [0,\infty)$.
\end{Fact}
\begin{Fact}\label{s3_F4}$\mbox{ \ }$
\begin{enumerate}
\item[{\bf (1)}] {\rm (Atiyah-Patodi-Singer \cite{APS})} 
There is an isomorphism: 
$$
{\mathcal H}_{\bar{g}}^2(X)\cong \Im \left[H^2_c(X;\R)\to H^2(X;\R)\right], 
$$
where $H^2_c(X;\R)$ denotes the second cohomology with compact support. 
\item[{\bf (2)}] {\rm (Dodziuk \cite{Do1})} 
For any $L^2_{\bar{g}}$-harmonic 2-form 
$\zeta\in {\mathcal H}_{\bar{g}}^2(X)$, 
there exists a sequence of closed $2$-forms 
$\{\zeta_j\} \subset \Omega_c^2(X)$ such that 
$$
\{
\begin{array}{l}
[\zeta_j]=[\zeta] \in H^2(X;\R), \\ \zeta_j\to \zeta \ \ \mbox{in the
$L^{k,2}_{\bar{g}}$-topology on $X$ for all $k\geq 0$ as $j\to
\infty$.}
\end{array}
\right.
$$
\end{enumerate} 
\end{Fact}
Recall that the Sobolev embedding 
$ L^{3,2}_{\bar{g}}(X) \subset C^{0,\alpha}(X)$ 
holds for all $0\leq \alpha <1$. 
This implies that 
$\sup_X|\zeta_j-\zeta|_{\bar{g}}\to 0$ as $j\to \infty$. 
Hence each closed $2$-form $\zeta_j$ is \emph{not harmonic} 
but \emph{almost-harmonic} for $j>>1$ unless $\zeta\equiv 0$.
\vspace{1mm}

\noindent
{\bf \ref{s3}.3. $L^2$-harmonic spinors.} 
>From now on we assume that 
$(X, \bar{g})$ is an oriented cylindrical $4$-manifold. 
An element $a\in H^2(X;\Z)$ is a \emph{characteristic element} 
if $a\equiv w_2(X)$ mod 2 
(where $w_2(X)\in H^2(X;\Z_2)$ is the second Stiefel-Whitney class).  
We will identify $a$ with its image in $H^2(X;\R)$.  
Let $P_Z : Z \times [1,\infty)\to Z\times \{1\}$ 
denote the canonical projection. 
We choose a characteristic element 
$a \in H^2(X; \Z) \cap \Im \left[H^2_c(X; \R) \to H^2(X; \R)\right]$ 
which is not a torsion class. 
Then there exists a Hermitian cylindrical $\C$-line bundle $\Ell$ over $X$ 
such that 
$$
\{
\begin{array}{l}
c_1(\Ell) = a,
\\
\Ell|_{Z\times [1,\infty)} = P^*_Z(\Ell|_{Z\times\{1\}}).
\end{array}
\right.
$$
Here $c_1(\Ell)$ is the first Chern class of $\Ell$. 

Let $\Sp^{\pm}_{\C}(a):= \Sp^{\pm}_{\C}\otimes \Ell^{1/2}$ 
denote the plus/minus spin bundle associated to $\Ell$. 
Here $\Sp^{\pm}_{\C}$ is the (virtual) plus/minus bundle. 
Then we have, as a Hermitian vector bundle, 
$$
\Sp^{\pm}_{\C}(a)|_{Z\times [1,\infty)} = 
P_Z^*(\Sp^{\pm}_{\C}(a)|_{Z\times \{1\}}).
$$
>From \cite[Theorem 2.7]{Do1}, there exists a unique 
$L^2_{\bar{g}}$-harmonic 2-form $\zeta\in {\mathcal H}_{\bar{g}}^2(X)$ 
such that its cohomology class $[\zeta]=a\in H^2(X;\R)$. 
In particular, from Fact \ref{s3_F4}-(2), 
there exists a sequence of closed $2$-forms 
$\{\zeta_j\} \subset \Omega_{\cp}^2(X)$ such that 
$$
\{
\begin{array}{l}
[\zeta_j] = a \in H^2(X;\R),
\\
\zeta_j\to \zeta \ \ \mbox{in} \ \ L^{3,2}_{\bar{g}}(X)\cap C^{0,\alpha}(X) 
\ \ 
\mbox{as} \ \ j\to\infty.
\end{array}
\right.
$$
Let ${\mathcal A}(\Ell)$ denote the space of $U(1)$-connections on $\Ell$. 
Then there exists a sequence of $U(1)$-connections 
$\{A_j\} \subset {\mathcal A}(\Ell)$ with compact support such that
$$ 
\frac{\sqrt{-1}}{2\pi} F_{A_j}:= \frac{\sqrt{-1}}{2\pi}dA_j = \zeta_j. 
$$ 
We also denote by $\dirac_{A_j} = \dirac_{A_j}^+:
\Gamma(\Sp^+_{\C}(a))\to \Gamma(\Sp^-_{\C}(a))$ the associated
(twisted) Dirac operator.
\vspace{2mm}

\noindent
{\bf Assumption A1.} 
>From now on we assume that, 
for a fixed cylindrical metric $\bar{g}_0\in \Riem^{\cyl}(X)$ with 
$\p_{\infty}\bar{g}_0=h$, 
$$
L^2\hbox{-}\ind \dirac_{A_j} : = 
\dim_{\C}(L^2\hbox{-}\Ker \dirac_{A_j}) - 
\dim_{\C}(L^2_{\ext}\hbox{-}\Ker \dirac_{A_j}^{\ *})> 0.
$$
Here $L^2\hbox{-}\ind \dirac_{A_j}$ stands for the $L^2$-index of 
$\dirac_{A_j}$ and $L^2_{\ext}\hbox{-}\Ker \dirac_{A_j}^{\ *}$ 
for the extended $L^2$-kernel of the adjoint operator $\dirac_{A_j}^{\ *}$ 
(cf.~\cite{APS, BW}). 
We emphasize that the above index is independent of $j$ 
(cf.~\cite{BW, N2}). 
\begin{Proposition}\label{s3_P2} 
For any $\bar{g}\in \Riem^{\cyl}(X)$ with $\p_{\infty}\bar{g}=h$, 
let $a \in H^2(X; \Z)$ be a non-torsion characteristic element and 
$\zeta\in {\mathcal H}_{\bar{g}}^2(X)$ 
its $L^2_{\bar{g}}$-harmonic representative. 
Then, under the assumption {\bf A1}, 
$$
Y^{\cyl}_{[\bar{g}]}(X)\leq \min
\{
4\pi\sqrt{2(a^+)^2}, \ Y^{\cyl}_{[h+dt^2]}(Z \times \R)
\}.
$$
Here $(a^+)^2:= \int_X|\zeta^+|_{\bar{g}}^2 d\sigma_{\bar{g}}$, and
$\zeta^+\in {\mathcal H}^{2,+}_{\bar{g}}(X)$ stands for 
the self-dual part of $\zeta$ with respect to $\bar{g}$. 
\end{Proposition}
\begin{Proof}
Recall that 
$Y^{\cyl}_{[\bar{g}]}(X)\leq Y^{\cyl}_{[h+dt^2]}(Z\times\R)$ 
(see \cite{AB2}). 
Hence it is enough to show that 
$Y^{\cyl}_{[\bar{g}]}(X)\leq 4\pi\sqrt{2(a^+)^2}$. 
Set $\omega:=4\sqrt{2}\zeta^+_j$ in Proposition \ref{s3_P1}. 
As we have seen, the form $\omega$ is not closed in general. 
We postpone the proof of the following assertion to the end of this section. 
\begin{Proposition}\label{s3_C1} 
For any $\bar{g}\in \Riem^{\cyl}(X)$ with $\p_{\infty}\bar{g}=h$, 
let $a \in H^2(X; \Z)$ be a non-torsion characteristic element 
and $\zeta \in {\mathcal H}_{\bar{g}}^2(X)$ 
its $L^2_{\bar{g}}$-harmonic representative. 
Then, under the assumption {\bf A1}, the following inequality holds: 
\begin{equation}\label{s3_eqq1} 
Y^{\cyl}_{[\bar{g}]}(X) \leq\|4\sqrt{2}\pi\zeta_j^+\|_{L_{\bar{g}}^2} 
= 4\pi\left(2 \int_X |\zeta_j^+|^2_{\bar{g}} d\sigma_{\bar{g}} 
\right)^{\frac{1}{2}}\quad \textrm{for any}\ j. 
\end{equation}
\end{Proposition}
By taking $j\to\infty$ in (\ref{s3_eqq1}), Proposition \ref{s3_C1}
implies that
$$
Y^{\cyl}_{[\bar{g}]}(X) \leq 4\pi\left(
2 \int_X |\zeta^+|^2_{\bar{g}} d\sigma_{\bar{g}}
\right)^{\frac{1}{2}} = 4\pi\sqrt{2(a^+)^2}.
$$
This completes the proof of Proposition \ref{s3_P2}.
\end{Proof}

\noindent
{\bf Assumption A2.} 
We assume that $b_2^-(X) = 0$. 
This implies that 
${\mathcal H}_{\bar{g}}^{2,+}(X) = {\mathcal H}_{\bar{g}}^2(X)$ 
for any $\bar{g}\in \Riem^{\cyl}(X)$. 
\begin{Theorem}\label{s3_T1}
Let $(X,\bar{g})$ be a cylindrical $4$-manifold modeled by $(Z,h)$. 
Then, under the assumptions {\bf A1} and {\bf A2},
$$
Y^{h\hbox{-}\cyl}(X)\leq \min\{4\pi\sqrt{2a^2},\ 
Y^{\cyl}_{[h+dt^2]}(Z\times \R)\}.
$$
\end{Theorem} 
\begin{lproof}{Proof}
Let $\bar{g}\in\Riem^{\cyl}(X)$ be any cylindrical metric 
with $\p_{\infty}\bar{g}=h$.  
For the 2-forms $\zeta$, $\zeta_j$ and 
the $U(1)$-connections $A_j$ as above, 
Fact \ref{s3_F5} and Assumption A1 imply that 
$L^2\hbox{-}\ind \dirac_{A_j,\bar{g}}>0$ for any $j$. 
Then, Proposition \ref{s3_P2} and Assumption A2 give the inequality 
$
Y^{\cyl}_{[\bar{g}]}(X) \leq 4\pi\sqrt{2(a^+)^2} = 4\pi\sqrt{2 a^2}.
$
Therefore,  $Y^{h\hbox{-}\cyl}(X) \leq 4\pi\sqrt{2 a^2}$. 
\end{lproof} 
Now we consider the following subspace of $\Riem(Z)$ 
$$
\Riem^*(Z):=\{h\in \Riem(Z) \ | \ \lambda({\mathcal L}_h)>0
\}\ (\ \subset \Riem^+(Z)\ ). 
$$
\begin{Fact}\label{s3_F5}{\rm (cf.~\cite{BW,N2})}
Let $A\in {\mathcal A}(\Ell)$ be a $U(1)$-connection 
with compact support $\supp(A) \subset X$. 
Let $h\simeq h^{\prime}$ be two metrics homotopic in $\Riem^*(Z)$. 
For any two cylindrical metrics 
$\bar{g}, \ \bar{g}^{\prime}\in \Riem^{\cyl}(X)$ 
with $\p_{\infty}\bar{g}=h$ and $\p_{\infty}\bar{g}^{\prime} = h^{\prime}$, 
then 
$L^2\hbox{-}\ind \dirac_{A,\bar{g}} = 
L^2\hbox{-}\ind \dirac_{A,\bar{g}^{\prime}}.$
\end{Fact} 
Then, Theorem \ref{s3_T1} combined with Fact \ref{s3_F5} 
implies the following. 
\begin{Corollary}\label{s3_Co1}
Let $h^{\prime}\in \Riem^*(Z)$ be a metric homotopic to $h$ in
$\Riem^*(Z)$. Then 
$$
Y^{h^{\prime}\hbox{-}\cyl} (X) \leq
\min
\{
4\pi\sqrt{2 a^2}, \ Y^{\cyl}_{[h^{\prime}+dt^2]}(Z^3\times \R)
\}.
$$
\end{Corollary} 
{\bf \ref{s3}.4. 
Behavior of $L^2$-harmonic spinors on the cylindrical ends.} 
Here we modify the Bochner techinque given in 
\cite{Hijazi}, \cite{Hitchin} and \cite{Lo} to our case. 
Let $\bar{g} \in \Riem^{cyl}(X)$ be a cylindrical metric 
with $\partial_{\infty}\bar{g} = h$ and 
$\check{g}= u^2 \cdot \bar{g}$ a conformal metric on $X$ 
with $u \in C_+^{\infty}(X)\cap L^{1,2}_{\bar{g}}(X)$. 
Then the Dirac operators 
$\dirac_{A,\check{g}}$ and $\dirac_{A,\bar{g}}$ on $\Sp^+_{\C}(a)$ 
are related as follows (cf.~\cite[Proposition 2]{Lo}): 
$$
\dirac_{A,\check{g}}(\psi) = u^{-\frac{5}{2}}\cdot \dirac_{A,\bar{g}}(
u^{\frac{3}{2}}\psi).
$$
For any $\psi\in \Ker\dirac_{A,\bar{g}}$, 
then we also obtain the harmonic spinor 
$\check{\psi}:= u^{-\frac{3}{2}}\cdot \psi \in \Ker\dirac_{A,\check{g}}$.  
Now let $\psi \in L^2_{\bar{g}} \hbox{-}\Ker\dirac_{A,\bar{g}}$ 
be a non-zero $L^2_{\bar{g}}$-harmonic plus-spinor.  
Then $\check{\psi}:= 
u^{-\frac{3}{2}} \cdot \psi\in \Ker\dirac_{A,\check{g}}$ 
is a non-zero harmonic spinor with respect to $(A, \check{g})$, 
and hence the Bochner formula gives that 
$$
0 = \dirac_{A,\check{g}}^{\ 2}\check{\psi} =-\Delta_{A,\check{g}}\check{\psi}
+\frac{1}{4}R_{\check{g}}\check{\psi} +\frac{1}{2}F_A^+\check{\ccd}
\check{\psi}.
$$
Here and below ``$\ccd$'' and ``$\check{\ccd}$'' stand for 
the Clifford multiplication corresponding to 
the metrics $\bar{g}$ and $\check{g}$, respectively.  
Then we have on $X(\ell) = X \setminus (Z \times (\ell,\infty))$ 
(for $\ell > 1$):
\begin{equation}\label{s3_eq1}
\begin{array}{rcl}
0 \!\!&=& \!\!\!\!
\displaystyle
\int_{X(\ell)}\left[
\<-\Delta_{A,\check{g}}\check{\psi},\check{\psi}\>+
\frac{1}{4}R_{\check{g}}|\check{\psi}|^2 +
\frac{1}{2}\<F_A^+\check{\ccd} \check{\psi}, \check{\psi}\>
\right]d\sigma_{\check{g}}
\\
\\
\!\!&=&\!\!\!\! 
\displaystyle
\int_{X(\ell)}\left[
|\nabla^{A,\check{g}}\check{\psi}|^2 +
\frac{1}{4}R_{\check{g}}|\check{\psi}|^2 +
\frac{1}{2}\<F_A^+\check{\ccd} \check{\psi}, \check{\psi}\>
\right]d\sigma_{\check{g}} 
\\
\\ 
&& \ \ \ \ \ \ \ \ \ \ \ \ \ \ \ \ \ \ \ \ \ \ \ \ \ \ \ 
\ \ \ \ 
\displaystyle
-
\int_{\p X(\ell)}\< \nabla^{A,\check{g}}_{\check{\nu}} \check{\psi},
\check{\psi}\> 
d\sigma_{\check{g}}|_{\p X(\ell)}
\\
\\
\!\!&\geq& \!\!\!\! 
\displaystyle
\int_{X(\ell)}\left[
|\nabla^{A,\check{g}}\check{\psi}|^2 +
\frac{1}{4} \left( R_{\check{g}} -2\sqrt{2} |F_A^+|_{\check{g}} 
\right)|\check{\psi}|^2 \right]d\sigma_{\check{g}} 
\\
\\
&& \ \ \ \ \ \ \ \ \ \ \ \ \ \ \ \ \ \ \ \ \ \ \ \ \ \ \ 
\ \ \ \ 
\displaystyle
-
\int_{\p X(\ell)}| \nabla^{A,\check{g}}_{\check{\nu}} \check{\psi}|
\cdot| \check{\psi}|
d\sigma_{\check{g}}|_{\p X(\ell)} \ .
\end{array}\!\!\!\!\!\!\!\!\!\!\!\!\!\!\!\!
\end{equation}
Here $\check{\nu}=u^{-1}\frac{\p}{\p t}$ is the outer unit normal 
vector field on $\p X(\ell)$ with respect to $\check{g}$.  
Since $\check{g}= u^2\cdot \bar{g}$, we have
$$
d\sigma_{\check{g}} = u^4 d\sigma_{\bar{g}} , \ \ \ \ 
|\check{\psi}| = u^{-\frac{3}{2}}|\psi|.
$$
Then, for certain positive constants $K$, $K^{\prime}$, 
we have the estimate:
\begin{equation}\label{s3_eq2}
\!\!
\begin{array}{rcl}
| \nabla^{A,\check{g}}_{\check{\nu}} \check{\psi}| \!\!&=& \!\!\!
u^{-1}\left| \nabla^{A,\check{g}}_{\p_t} (u^{-\frac{3}{2}}\psi)\right|
= u^{-1}\left| -\frac{3}{2}(u^{-\frac{5}{2}}\p_t u)\psi + u^{-\frac{3}{2}} 
\nabla^{A,\check{g}}_{\p_t}{\psi}\right|
\\
\\
\!\!&=& \!\!\!
u^{-\frac{5}{2}}
\left|
-\frac{3}{2}( u^{-1}\p_tu)\psi + \nabla^{A,\bar{g}}_{\p_t}{\psi}
-\frac{1}{4}\sum_{j=1}^4 e_j(u)\left[(\frac{\p}{\p t})\ccd,e_j\ccd\right]
\psi\right|
\\
\\
&\leq& \!\!\!
K \cdot u^{-\frac{5}{2}}
\left(
 u^{-1}|\p_tu|\cdot |\psi| + |\nabla^{A,\bar{g}}_{\p_t}{\psi}|
+ u^{-1} |du|_{\bar{g}}\cdot  |\psi|
\right)
\\
\\
&\leq&\!\!\!
K^{\prime} \cdot u^{-\frac{5}{2}}
\left(
 u^{-1}|du|_{\bar{g}}\cdot  |\psi| + |\nabla^{A,\bar{g}}_{\p_t}{\psi}|
\right) ,
\end{array}\!\!\!\!\!\!\!\!\!\!\!\!\!\!\!\!\!\!\!\!\!\!\!\!
\end{equation}
where $\{e_1=\frac{\p}{\p t}, e_2,e_3,e_4\}$ is a local orthonormal frame 
on $\p X(\ell)$ with respect to $\bar{g}$.  
We use (\ref{s3_eq2}) to obtain that 
\begin{equation}\label{s3_eq3}
\!\!\begin{array}{l}
\displaystyle
\int_{\p X(\ell)}\!\!| \nabla^{A,\check{g}}_{\check{\nu}} \check{\psi}|
\cdot| \check{\psi}|
d\sigma_{\check{g}}|_{\p X(\ell)} 
\\
\\
\displaystyle
\quad \leq
K \int_{\p X(\ell)} u^{-\frac{5}{2}}\left(
 u^{-1}|du|_{\bar{g}}\cdot  |\psi| + |\nabla^{A,\bar{g}}_{\p_t}{\psi}|
\right)u^{-\frac{3}{2}} |\psi| u^4 d\sigma_{\bar{g}}|_{\p X(\ell)}
\\
\\
\quad \leq \displaystyle
K \int_{\p X(\ell)}\!\! \left(
u^{-1}|du|_{\bar{g}}\cdot  |\psi|^2 + |\nabla^{A,\bar{g}}_{\p_t}{\psi}|
\cdot  |\psi|
\right)d\sigma_{\bar{g}}|_{\p X(\ell)}\ .
\end{array} 
\end{equation}
\begin{Lemma}\label{s3_L5}
Let $u\in C_+^{\infty}(X)\cap L_{\bar{g}}^{1,2}(X)$ be a function 
satisfying $\L_{(\bar{g},\omega)}u= \lambda_{(\bar{g},\omega)} u$ 
on $Z\times [\ell_0,\infty)$ for a fixed $\ell_0\geq 1$. 
Then  
$$
B_{\ell}:=\int_{\p X(\ell)} \left(
u^{-1}|du|_{\bar{g}}\cdot  |\psi|^2 + |\nabla^{A,\bar{g}}_{\p_t}{\psi}|
\cdot  |\psi|
\right)d\sigma_{\bar{g}}|_{\p X(\ell)}\to 0 \ \ \mbox{as $\ell\to \infty$}.
$$
\end{Lemma}
\begin{lproof}{Proof of Proposition \ref{s3_C1}}
Set $\omega:=4\sqrt{2}\pi\zeta^+_j = \sqrt{-1}\cdot 2\sqrt{2}
F_{A_j}^+$ and let $u$ be the same function as in Lemma \ref{s3_L4}.
Suppose that the $C^{2,\alpha}$-metric $\check{g}=u^2\cdot \bar{g}$
satisfies that $R_{\check{g}}> |\omega|_{\check{g}}$. Perturbing $u$
on a compact set, we may assume that $u$ is smooth and that 
$R_{\check{g}} > |\omega|_{\check{g}}$.
Then, the inequality (\ref{s3_eq1}) combined with 
the estimate (\ref{s3_eq3}) and Lemma \ref{s3_L5} implies that 
\begin{equation}\label{s3_eq4}
0\geq \int_X \left[|\nabla^{A_j,\check{g}}\check{\psi}|^2 +
\frac{1}{4} \left( R_{\check{g}} -2\sqrt{2} |F_{A_j}^+|_{\check{g}} 
\right)|\check{\psi}|^2 \right]d\sigma_{\check{g}} \ .
\end{equation}
Hence the condition 
$R_{\check{g}}>|\omega|_{\check{g}} = 2\sqrt{2}|F_{A_j}^+|_{\check{g}}$ 
contradicts the inequality (\ref{s3_eq4}). 
Finally, from Proposition \ref{s3_P1},
$Y^{\cyl}_{[\bar{g}]}(X)\leq \|2\sqrt{2}F_{A_j}^+\|_{L_{\bar{g}}^2}=
\|4\sqrt{2}\pi\zeta_j^+\|_{L_{\bar{g}}^2}$.
\end{lproof}
Now we have to prove Lemma \ref{s3_L5}. 
There are two independent tasks here to deal with, i.e., 
the decay of the function $u$ and the pointwise norms $|\psi|$, 
$|\nabla^{A,\bar{g}}_{\p_t}{\psi}|$. 
The spinor norms are taken care by the following result. 
\begin{Fact}\label{s3_F6}{\rm (\cite{APS})}
Let $\psi$ be an $L_{\bar{g}}^2$-harmonic spinor as above. 
Then the pointwise norms $|\psi|$, $|\nabla^{A,\bar{g}}_{\p_t}{\psi}|$ 
decay exponentially on the cylinder $Z\times [1,\infty)$. 
\end{Fact}
To complete the proof of Lemma \ref{s3_L5}, 
we need the following technical result. 
\begin{Lemma}\label{s3_L6} 
There exists a positive constant $K$ such that 
$u^{-1}|du|_{\bar{g}}\leq K$ on the cylinder $Z\times [1,\infty)$. 
\end{Lemma}
{\em Proof of Lemma \ref{s3_L5}.} 
By Fact \ref{s3_F6} and Lemma \ref{s3_L6}, 
there exist positive constants $K_1$, $K_2$, $K$ and $\kappa$ such that 
$$
\begin{array}{cr}
B_{\ell}\leq K_1 \!\cdot\!\Vol_h(Z)\left(\!
K_2\cdot |\psi|^2+ |\nabla^{A,\bar{g}}_{\p_t}{\psi}|\cdot|\psi|\!
\right)\! \leq\! K\!\cdot 
\!e^{-\kappa\ell} \to 0 
& 
\end{array}
$$
as $\ell\to \infty$.\hfill $\Box$
%
\begin{lproof}{Proof of Lemma \ref{s3_L6}}
Recall that the support $\supp(\omega)\ (\subset X)$ is compact. 
Hence there exists $\ell_1$ ($\geq \ell_0$) such that 
the restriction of $\omega$ is zero 
on the cylinder $Z\times [\ell_1-1,\infty)$. 
In particular, the operator 
$\L_{(\bar{g},\omega)}= \L_{\bar{g}}$ on $Z\times [\ell_1-1,\infty)$. 
Recall that $\lambda_h$ and $\lambda_{(\bar{g}, \omega)}$ 
denote respectively the first eigenvalue of the operator ${\mathcal L}_h$ 
and the bottom of the spectrum of $\L_{(\bar{g},\omega)}$. 
Then, we consider the following operator $L_h$ on $Z$ 
$$
L_h:= -6\Delta_h + (R_h - \lambda_{(\bar{g}, \omega)}) = 
{\mathcal L}_h - \lambda_{(\bar{g}, \omega)}. 
$$
Let $\{(\mu_j,\phi_j)\}$ be the eigenvalues and eigenfunctions of $L_h$, 
i.e., 
$$
L_h\phi_j= \mu_j\phi_j \ \ \ \mbox{with}\ \ \ 
\int_Z \phi_j \cdot \phi_k d\sigma_h = \delta_{jk}
\ \ \ \mbox{for} \ \ \  j,k=1,2,\ldots, \ \  \ 
$$
and $\mu_1<\mu_2\leq\mu_3\leq\cdots$. 
Here we may assume that $\phi_1>0$ on $Z$. 
Since $\lambda_h > \lambda_{(\bar{g}, \omega)}$, 
we obtain that $\mu_1 > 0$. 

Let $(z,t)\in Z\times [\ell_1,\infty)$ 
be a product coordinate system associated to $\bar{g}$. 
Note that 
$\L_{(\bar{g},\omega)}= 
-6\cdot\p_t^2 + {\mathcal L}_h$ on $Z \times [\ell_1,\infty)$. 
Then on the cylinder $Z \times [\ell_1,\infty)$, 
$$
\!\!\!\begin{array}{l} 
(\L_{(\bar{g},\omega)}\! -\!\lambda_{(\bar{g}, \omega)}) 
(e^{-\sqrt{\mu_j/6}\cdot t}\!\phi_j(z))
= 
L_h(e^{-\sqrt{\mu_j/6}\cdot t}\!\phi_j(z))\! - \!
6\!\cdot\!\p_t^2(e^{-\sqrt{\mu_j/6} \cdot t}\!\phi_j(z)) 
\\
\\
\ \ \ \ \ \ \ \ \ \ \ =
\displaystyle
\mu_j\cdot (e^{-\sqrt{\mu_j/6} \cdot t}\phi_j(z)) - 
\mu_j\cdot (e^{-\sqrt{\mu_j/6} \cdot t}\phi_j(z)) = 0.
\end{array}
$$
Hence $\L_{(\bar{g},\omega)}(e^{-\sqrt{\mu_j/6}\ \cdot t}\phi_j(z)) = 
\lambda_{(\bar{g},\omega)} \cdot e^{-\sqrt{\mu_j/6}\ \cdot t}\phi_j(z)$ 
on $Z \times [\ell_1,\infty)$. 
Then the function $u\in C^{\infty}_+(X)\cap L_{\bar{g}}^{1,2}(X)$ 
restricted on $Z\times [\ell_1,\infty)$ can be represented as 
$$
\{
\begin{array}{l}
u(z,t)= \sum_{j\geq 1} a_j \cdot e^{-\sqrt{\mu_j/6}\ \cdot t}\phi_j(z),
\\
\sum_{j\geq 1} a_j^2 < \infty.
\end{array}
\right.
$$
We set $\mu:=\mu_{j_0}= \min\{ \mu_j \ | \ a_j \neq 0\}$ and
$\upsilon_j:= \sqrt{\mu_j/6}-\sqrt{\mu/6}>0$ for $\mu_j>\mu$.  Let
$k+1$ ($\geq 1$) be the multiplicity of the $j_0$-th eigenvalue.  Then
$$
\begin{array}{l}
\displaystyle
\!\!
u(z,t)\!= \!e^{-\sqrt{\mu/6} \cdot t}
\left(\!
a_{j_0} \phi_{j_0}(z)\!+\!\cdots\!+ \!a_{j_0+k} \phi_{j_0+k}(z)\!
+\!\sum_{\mu_j>\mu}\!\!a_j \!\cdot \!e^{-\upsilon_j\cdot t}\!\phi_j(z)\!
\right).
\end{array}
$$
It then follows that there exists a positive constant $K$
such that 
$$
|\p_t u| \leq K\cdot e^{-\sqrt{\mu/6}\ \cdot t}
, \ \ \ \ \ 
|\nabla^Zu|_h \leq K \cdot e^{-\sqrt{\mu/6}\ \cdot t} \ .
$$
This implies the estimate
\begin{equation}\label{s3_eq5}
|du|_{\bar{g}}\leq \bar{K}\cdot e^{-\sqrt{\mu/6}\ \cdot t}
\end{equation}
for a constant $\bar{K}>0$. 
Now we set 
$v(z):= a_{j_0} \phi_{j_0}(z)+\cdots+ a_{j_0+k} \phi_{j_0+k}(z)$. 
Here we need the following fact. 
\begin{Claim}\label{s3_CS}
With the above understood, then 
$$
\{
\begin{array}{l}
\mu=\mu_1 \ \ \mbox{is the first eigenvalue of $L_h$},
\\
v(z)= a_1\phi_1(z) > 0 \ \ \mbox{on $Z$}.
\end{array}
\right.
$$
\end{Claim}
\begin{Proof}
Note that if there exists $z_0\in Z$ such that $v(z_0)<0$, 
then $u(z_0,t)<0$ for sufficiently large $t>>\ell_1$. 
This contradicts the positivity of $u$ everywhere on $X$. 
Hence $v(z)\geq 0$ and $v(z)\not\equiv 0$ on $Z$. 
Now suppose that $\mu$ is not the first eigenvalue.  
Note that the function $v(z)$ itself is 
an eigenfunction corresponding to $\mu$. 
The condition $\mu \not= \mu_1$ implies that 
$\int_Z \phi_1(z) \cdot v(z)d\sigma_h =0$, 
which contradicts that $\phi_1(z)>0$, $v(z)\geq 0$ and 
$v(z)\not\equiv 0$ on $Z$.  
Therefore, $\mu=\mu_1$.
\end{Proof} 
Note that Claim \ref{s3_CS} implies the estimate 
\begin{equation}\label{s3_eq6} 
0<\tilde{K}^{-1} \cdot e^{-\sqrt{\mu/6}\ \cdot t} \leq u(z,t) \leq 
\tilde{K} \cdot e^{-\sqrt{\mu/6}\ \cdot t} 
\end{equation}
on $Z\times [\ell_1,\infty)$ for a constant $\tilde{K}>0$. 
Then the estimates (\ref{s3_eq5}) and (\ref{s3_eq6}) imply that 
$u^{-1}|du|_{\bar{g}}\leq K$ on $X$ for a constant $K > 0$.
This completes the proof of Lemma \ref{s3_L6}.
\end{lproof}

\section{Examples}\label{s4} 
Let $\C\P^2$ denote the complex projective plane. 
In \cite{Le} (cf.~\cite{GL}), LeBrun proved that 
$Y(\C\P^2) = 12\sqrt{2}\pi (< 8\sqrt{6}\pi = Y(S^4))$, 
and that $Y(\C\P^2)$ is achieved by a conformal class $C$ 
if and only if $C$ coincides with the pullback $\Phi^*[g_{FS}]$ 
of the conformal class $[g_{FS}]$ by a diffeomorphism $\Phi$. 
Here $g_{FS}$ denotes the Fubini-Study metric on $\C\P^2$. 
In this section, we give some estimates of $Y^{h\textrm{-}\cyl}(X)$ 
for particular cylindrical $4$-manifolds $X$, 
which are generalizations of part of the above result. 

\noindent
{\bf \ref{s4}.1. Example~1.} 
Let denote by $\Ell\to \C\P^1$ the anti-canonical $\C$-line bundle 
over the complex projective line $\C\P^1$ 
and by $\Ell_{\ell}:=\Ell^{\otimes\ell}\to \C\P^1$ 
the $\C$-line bundle of degree $\ell\geq 1$. 
Then we choose a Hermitian metric $H$ on the bundle $\Ell_{\ell}$, 
and consider the disk and the circle bundles of $\Ell_{\ell}$: 
$$ 
\begin{array}{rcccl}
D_{\ell}&:=& D(\Ell_{\ell})&=& \{ v\in \Ell_{\ell} \ | \ H(v,v)\leq 1 \}, 
\\
Z_{\ell}&:=& S(\Ell_{\ell})&=& \p D(\Ell_{\ell}) \cong S^3/\Gamma_{\ell} 
=  L(\ell,1).
\end{array} 
$$ 
Here $\Gamma_{\ell}= \Z/\ell\Z\subset \{ \xi\in \C \ | \ |\xi|=1 \}$ 
with the generator $ \xi_{\ell}= e^{\frac{2\pi\sqrt{-1}}{\ell}}$ 
and the group $\Gamma_{\ell}$ acts on the sphere 
$S^3 = S^3(1) \subset \C^2$ by 
$(w^1,w^2) \mapsto (\xi_{\ell}\cdot w^1,\xi_{\ell}\cdot w^2)$. 
We also denote by $X_{\ell}$ the total space of $\Ell_{\ell}$. 
We consider $X_{\ell}$ as an open $4$-manifold with a tame end: 
$$
X_{\ell} \cong D_{\ell}\cup_{Z_{\ell}} (Z_{\ell}\times [0,\infty)). 
$$
First we note some necessary facts.
\begin{Fact}\label{s4_F1}$\mbox{ \ }$ 
\begin{enumerate}
\item[{\bf (1)}] 
When $\ell=1$, then $X_1$ is diffeomorphic to $\C\P^2 \setminus \{q\}$, 
and hence $X_1$ is an open $4$-manifold with tame end 
$S^3 \times [0, \infty)$, 
where $q$ is a point in $\C\P^2$. 
When $\ell \geq 2$, then the one-point compactification 
$M_{\ell} := X_{\ell} \cup \{\check{p}_{\infty}\}$ 
has a natural orbifold structure 
with $\Sigma_{\Gamma} = \{(\check{p}_{\infty}, \Gamma_{\ell})\}$. 
\item[{\bf (2)}] 
{\rm (cf.~\ \cite[Example 4.1.27]{N2})}\qquad  
There exists a cylindrical metric $\bar{g}_0 \in $ 
$\Riem^{\cyl}(X_{\ell})$ with 
$\p_{\infty}\bar{g}_0 = h_{\Gamma_{\ell}} \in \Riem(Z_{\ell})$ such that 
$R_{\bar{g}_0} > 0$ on $X_{\ell}$. 
\item[{\bf (3)}] 
$H^2(X_{\ell};\Z)\cap \Im \left[H_c^2(X_{\ell};\R) \to 
H^2(X_{\ell};\R)\right] \cong \Z$, 
and $b_2^-(X_{\ell}) = 0$ for an appropriate orientation of $X_{\ell}$. 
\end{enumerate}
\end{Fact} 
>From \cite[Lemma~2.14]{AB2}, we note that 
$Y^{h_S\textrm{-}cyl}(X_1) = Y(\C\P^2)$, 
where $h_S = h_{\Gamma_1} \in \Riem(S^3)$ denotes 
the standard metric of constant curvature one. 
With these understood, we prove the following more general result. 
\begin{Theorem}\label{s4_T1} 
\quad 
\begin{enumerate} 
\item[{\bf (1)}] 
$$ 
0< Y^{h\hbox{-}\cyl}(X_{\ell}) \leq 
\min\{4(\ell+2)\sqrt{\frac{2}{\ell}}\ \pi,\ 
Y^{\cyl}_{[h+dt^2]}(Z_{\ell} \times \R)\} 
$$ 
for $\ell\geq 1$ and any metric $h\in \Riem^*(Z_{\ell})$ 
homotopic to $h_{\Gamma_{\ell}}$ in $\Riem^*(Z_{\ell})$. 
\item[{\bf (2)}] 
For any metric $h \in \Riem(Z_{\ell})$ 
which is sufficiently $C^2$-close to $h_{\Gamma_{\ell}}$, 
$$
0< Y^{h\hbox{-}\cyl}(X_{\ell})\ \leq \ 
\{ 
\begin{array}{rll} 
12\sqrt{2}\pi = Y(\C\P^2) & \mbox{if} \ \ \ \ell=1, 
\\ 
Y^{\cyl}_{[h+dt^2]}(Z_{\ell} \times \R)\ \ 
& \mbox{if} \ \ \ \ell\geq 2. 
\end{array} 
\right. 
$$ 
In particular, 
$0 < Y^{\orb}(M_{\ell}) = 
Y^{h_{\Gamma_{\ell}}\hbox{-}\cyl}(X_{\ell}) \leq  
8\sqrt{\frac{6}{\ell}}\ \pi$ 
\ for $\ell \geq 2$. 
\end{enumerate} 
\end{Theorem} 
\begin{Remark} {\rm We note the following:} 
\begin{enumerate} 
\item[{\bf (1)}] {\rm The inequality $4(\ell+2)\sqrt{\frac{2}{\ell}}\ 
\pi<8\sqrt{6}\pi\ (\ = Y(S^4)\ )$ holds only if $1\leq \ell\leq 7$.}
\item[{\bf (2)}] {\rm Let $M$ be a closed $4$-manifold. 
For any $\ell\geq 2$ and a positive integer $k$, then} $
Y^{\orb}(M\#k M_{\ell}) \leq 8\sqrt{\frac{6}{\ell}}\ \pi.$ 
\item[{\bf (3)}] {\rm If we replace $X_1 = \C\P^2 \setminus \{q\}$ 
by $\C\P^2 \setminus \{q_1, \ldots, q_s\}$, 
a similar result still holds for $\C\P^2 \setminus \{q_1, \ldots, q_s\}$.} 
\item[{\bf (4)}] {\rm Unfortunately, the estimate 
$Y^{h\textrm{-}\cyl}(X) \leq 4\pi\sqrt{2a^2}$ 
in Theorem \ref{s3_T1} is not effective for 
$X=X_{\ell}$ ($\ell\geq 2$) and $h=h_{\Gamma_{\ell}}$ since 
$4(\ell+2)\sqrt{\frac{2}{\ell}}\pi > 8 \sqrt{\frac{6}{\ell}}\pi$ 
for $\ell \geq 2$.}
\end{enumerate}
\end{Remark}
\begin{lproof}{Proof of Theorem \ref{s4_T1}} 
For $\ell \geq 1$, let ${h}\in \Riem^*(Z_{\ell})$ be a metric 
which is homotopic to $h_{\Gamma_{\ell}}$ in $\Riem^*(Z_{\ell})$. 
By Fact \ref{s4_F1}-(2), it is easy to construct 
a cylindrical metric $\bar{g}$ on $X$ satisfying 
$R_{\bar{g}} > 0$ on 
$X_{\ell}(1) = X_{\ell} \setminus (Z_{\ell} \times (1, \infty))$, 
$\bar{g} = h + dt^2$ on $Z_{\ell} \times [2, \infty)$ and 
$[\bar{g}] = [h + dt^2]$ on $Z_{\ell} \times [1, \infty)$. 
Then, from \cite[Proposition~4.6, Theorem~4.7]{AB2}, 
we also obtain that 
$$ 
Y^{\cyl}_{[\bar{g}]}(X_{\ell}) \geq \min \{Y_1, Y_2\} > 0, 
$$ 
and hence $Y^{h\textrm{-}\cyl}(X_{\ell}) > 0$. 
Here $Y_1 = 
Y_{[\bar{g}|_{X_{\ell}(1)}]}
(X_{\ell}(1), \p X_{\ell}(1); [\bar{g}|_{\p X_{\ell}(1)}]) > 0$ 
denotes the relative Yamabe constant of 
$(X_{\ell}(1), \p X_{\ell}(1); [\bar{g}|_{X_{\ell}(1)}])$ (see~\cite{AB0}) 
and 
$$ 
\begin{array}{rcl} 
Y_2 &=& Y^{\cyl}_{[\bar{g}|_{Z_{\ell} \times [1, \infty)}]}
(Z_{\ell} \times [1, \infty), Z_{\ell} \times \{1\}; 
[\bar{g}|_{Z_{\ell} \times \{1\}}]) 
\\ 
\\ 
&=& Y^{\cyl}_{[h + dt^2]}(Z_{\ell} \times [1, \infty), 
Z_{\ell} \times \{1\}; [\bar{g}|_{Z_{\ell} \times \{1\}}])  > 0
\end{array}
$$ 
denotes the relative cylindrical Yamabe constant of 
$(Z_{\ell} \times [1, \infty), Z_{\ell} \times \{1\}; [h + dt^2])$ 
(see~\cite[Section~4.1]{AB2}), respectively. 
>From Theorem \ref{s3_T1}, Corollary \ref{s3_Co1} and Fact \ref{s4_F1}-(3), 
in order to prove the assertion {\bf (1)}, 
it is enough to show that there exists a characteristic element 
$a \in H^2(X_{\ell}; \Z) \cap 
\Im \left[H_c^2(X_{\ell}; \R) \to H^2(X_{\ell};\R)\right]$ 
satisfying the assumption {\bf A1} 
for the cylindrical metric $\bar{g}_0 \in \Riem^{\cyl}(X_{\ell})$ 
in Fact \ref{s4_F1}-(2) 
and $a^2 = \frac{(\ell + 2)^2}{\ell}$. 

Note that $X_{\ell}$ has a natural complex structure. 
Then let 
$$
a_0 \in H^2(X_{\ell}; \Z) \cap \Im \left[H_c^2(X_{\ell}; \R) \to
H^2(X_{\ell}; \R) \right] \cong \Z
$$ 
be the generator satisfying $c_1(X_{\ell})=(\ell+2)a_0$, 
and set $a:= (\ell+2)a_0$. 
It is easy to see that $a \equiv w_2(X_{\ell})$ mod $2$. 
Then, 
there exist the $L^2_{\bar{g}_0}$-harmonic $2$-form 
$\zeta \in {\mathcal H}^2_{\bar{g}_0}(X_{\ell}) 
= {\mathcal H}^{2,+}_{\bar{g}_0}(X_{\ell})$ 
and a sequence of closed $2$-forms 
$\{\zeta_j\} \subset \Omega^2_{\cp}(X_{\ell})$ 
such that 
$$
\{
\begin{array}{l}
[\zeta_j]=[\zeta] = a \in H^2(X_{\ell}; \Z), 
\\ 
\zeta_j\to \zeta \ \
\mbox{in $L^{k,2}_{\bar{g}_0}$-topology for any $k\geq 1$}.
\end{array} 
\right.
$$
Now we consider the determinant $\C$-line bundle $\det(a)$ 
and the associated plus/minus spin bundle 
$\Sp^{\pm}_{\C}(a)= \Sp^{\pm}_{\C}\otimes \det(a)^{1/2}$. 
Then we may assume that, as Hermitian vector bundles, 
$$
\det(a)|_{Z_{\ell}\times[1,\infty)}\!
=\! P^*_{Z_{\ell}}\!(\det(a)|_{Z_{\ell}\times\{1\}}\!), \ \ 
\Sp^{\pm}_{\C}\!(a)|_{Z_{\ell}\times[1,\infty)} \!
\!= P^*_{Z_{\ell}}(\Sp^{\pm}_{\C}(a)|_{Z_{\ell}\times\{ 1\} }\!),
$$
where $P_{Z_{\ell}}: Z_{\ell}\times[1,\infty)\to Z_{\ell}\times\{ 1\}$ 
denotes the canonical projection. 
Moreover, there exist $U(1)$-connections 
$A_j\in {\mathcal A}(\det(a))$ with compact support such that 
$\frac{\sqrt{-1}}{2\pi}F_{A_j} = \frac{\sqrt{-1}}{2\pi}dA_j = \zeta_j$. 
Now we denote by 
$$
\dirac_{A_j} := \dirac_{A_j,\bar{g}_0}^{\ +} : 
\Gamma(\Sp^+_{\C}(a)) \lra \Gamma(\Sp^-_{\C}(a))
$$
the corresponding twisted Dirac operator. 
Recall that 
$\p_{\infty}\bar{g}_0 = h_{\Gamma_{\ell}} \in \Riem(Z_{\ell})$, 
and hence 
$\lambda({\mathcal L}_{h_{\Gamma_{\ell}}}) > 0$. 
Then the Hirzebruch signature
formula (cf.~\cite[Example~4.1.9]{N2}) combined with
$\lambda({\mathcal L}_{h_{\Gamma_{\ell}}})>0$ gives
$$
\begin{array}{rcl}
L^2\hbox{-}\ind \ \dirac_{A_j} &=& \dim_{\C}(L^2\hbox{-}\Ker
\dirac_{A_j,\bar{g}_0}^{\ +}) - \dim_{\C}(L^2_{\ext}\hbox{-}\Ker
\dirac_{A_j,\bar{g}_0}^{\ -}) 
\\ 
\\ 
&=& 
\frac{1}{8}\int_{X_{\ell}}c_1(A_j)\wedge c_1(A_j) -\frac{1}{8}
\int_{X_{\ell}}\frac{1}{3} p_1(\nabla^{\bar{g}_0}) 
- \frac{1}{2}\eta(h_{\Gamma_{\ell}}).
\end{array}
$$
Here $\eta(h_{\Gamma_{\ell}})$ stands for the eta invariant 
of the Dirac operator
$$
\begin{array}{c}
\left(
\frac{\p}{\p t}
\right)\ccd \dirac_{(Z_{\ell},h_{\Gamma_{\ell}})}: 
\Gamma(\Sp^+_{\C}|_{Z_{\ell}\times \{1\}}) \lra 
\Gamma(\Sp^+_{\C}|_{Z_{\ell}\times \{1\}})
\end{array}
$$
and $\dirac_{(Z_{\ell},h_{\Gamma_{\ell}})}$ the spin Dirac operator on 
$(Z_{\ell}\cong Z_{\ell}\times \{1\},h_{\Gamma_{\ell}})$, 
respectively. 
\begin{Claim}\label{s4_C}
$L^2\hbox{-}\ind \ \dirac_{A_j}=1$. 
\end{Claim}
{\em Proof.}
Indeed, we have 
$$
\{
\begin{array}{l}
\int_{X_{\ell}}c_1(A_j)\wedge c_1(A_j) =
\frac{(\ell+2)^2}{\ell} \ \ \mbox{and}
\\
\\
\int_{X_{\ell}}\frac{1}{3} p_1(\nabla^{\bar{g}_0}) = \tau(X_{\ell}) + 
\eta_{\sign}(h_{\Gamma_{\ell}}) = 1+ \eta_{\sign}(h_{\Gamma_{\ell}}),
\end{array}
\right.
$$
where $\eta_{\sign}(h_{\Gamma_{\ell}})$ is the signature defect 
(see \cite{APS} or \cite[(4.1.34)]{N2}). 
Note that $4\cdot\eta(h_{\Gamma_{\ell}}) 
+ \eta_{\sign}(h_{\Gamma_{\ell}})= 
- \frac{4(\ell-1)}{\ell} + \ell-1$ 
(see \cite{N1} or \cite[Example~4.1.27]{N2}). 
Then we have
$$
\begin{array}{rlr}
L^2\hbox{-}\ind \ \dirac_{A_j} =& 
\frac{1}{8} 
\left[ 
\frac{(\ell+2)^2}{\ell} -1 - 
\left(
4\cdot\eta(h_{\Gamma_{\ell}}) + \eta_{\sign}(h_{\Gamma_{\ell}}) 
\right) 
\right]& 
\\ 
=& 
\frac{1}{8\ell} 
\left( 
\ell^2+4\ell+4-4\ell + 4 -\ell^2 
\right) = 1. & \ \ \Box 
\end{array} 
$$
We also note that 
$a^2=\int_{X_{\ell}}\zeta\wedge\zeta = \frac{(\ell+2)^2}{\ell}$ 
by the choice of $a$. This completes the proof of the assertion {\bf (1)}. 

By the Remark after \cite[Proposition~6.5]{AB2}, 
$$
\begin{array}{c}
Y^{\cyl}_{[h_{\Gamma_{\ell}} + dt^2]}(Z_{\ell} \times \R) 
= \frac{Y(S^4)}{|\Gamma_{\ell}|^{1/2}} = 8\sqrt{\frac{6}{\ell}}\pi.
\end{array}
$$
>From the continuity of Yamabe constants \cite[Proposition~3.4]{ABKS}, 
Remark~4.1-(4), and that 
$Y^{\cyl}_{[h_S + dt^2]}(S^3 \times \R) = Y(S^4) 
= 8\sqrt{6}\pi > 12\sqrt{2}\pi$, 
we obtain that 
$$
Y^{h\hbox{-}\cyl}(X_{\ell})\ \leq \ 
\{ 
\begin{array}{rll} 
12\sqrt{2}\pi = Y(\C\P^2) & \mbox{if} \ \ \ \ell=1, 
\\ 
Y^{\cyl}_{[h+dt^2]}(Z_{\ell} \times \R)\ \ & \mbox{if} \ \ \ \ell\geq 2 
\end{array} 
\right. 
$$ 
for any metric $h \in \Riem(Z_{\ell})$ 
which is sufficiently $C^2$-close to $h_{\Gamma_{\ell}}$. 
This completes the proof of the assertion {\bf (2)}. 
\end{lproof}

\begin{Corollary}\label{s4_C2}{\rm \!(cf.~\cite[Theorem~B]{GL})}\!
Set $X_{\ell,k,m} = X_{\ell} \# k\C\P^2 \# m(S^1\!\times\! S^3)$. 
For any integers $\ell\geq 1$ and $k, \ m \geq 0$, 
\begin{equation}\label{s4_eq1}
\!
\begin{array}{c}
0\!<  
\!Y^{h\hbox{-}\cyl}(X_{\ell,k,m})
\leq \min
\{\!4\pi\sqrt{2(k+\frac{(\ell+2)^2}{\ell})},\ 
Y^{\cyl}_{[h+dt^2]}(Z_{\ell}\times \R)
\!\}
\end{array}\!\!\!\!\!\!\!\!\!\!
\end{equation}
for any metric $h\in \Riem^*(Z_{\ell})$ homotopic to 
$h_{\Gamma_{\ell}}$ in $\Riem^*(Z_{\ell})$.
\end{Corollary}
\begin{Remark}
{\rm 
We note that 
$4\pi\sqrt{2(k+\frac{(\ell+2)^2}{\ell})} < Y(S^4)=8\sqrt{6}\pi$ 
\ \ if $k< 12 - \frac{(\ell+2)^2}{\ell} \leq 4$. 
} 
\end{Remark}
\begin{lproof}{Proof of Corollary \ref{s4_C2}}
There exists a metric $\tilde{g}\in \Riem^{\cyl}(X_{\ell,k,m})$ such that 
$R_{\tilde{g}}>0$ and $\p_{\infty}\tilde{g} = h_{\Gamma_{\ell}}$. 
Then, similarly to the case of $X_{\ell}$, 
we obtain that 
$Y^{h\hbox{-}\cyl}(X_{\ell,k,m})>0$ for any metric
$h\in\Riem^*(Z_{\ell})$ homotopic to 
$h_{\Gamma_{\ell}}$ in $\Riem^*(Z_{\ell})$. 
Note that 
$$
\begin{array}{l}
{\mathcal H}_{\tilde{g}}^{2,+}(X_{\ell,k,m})\cap 
H^2(X_{\ell,k,m};\Z) =
{\mathcal H}_{\tilde{g}}^{2}(X_{\ell,k,m})\cap 
H^2(X_{\ell,k,m};\Z) 
\\
\\
\ \ \ \ \ \ \  
\cong 
\left({\mathcal H}_{\tilde{g}}^{2}(X_{\ell})\cap 
H^2(X_{\ell};\Z)\right) \oplus H^2(\C\P^2;\Z)^{\oplus k} \cong \Z \oplus
\Z^{\oplus k}.
\end{array}
$$
Then we choose $\tilde{a}=(a,1,\cdots,1)\in \Z \oplus \Z^{\oplus k}$. 
Similarly to the case of $X_{\ell}$, 
we also obtain that $L^2\hbox{-}\ind\dirac_{\tilde{A}}=1$, 
where $\dirac_{\tilde{A}}$ denotes the corresponding 
Dirac operator on $\Sp^+_{\C}(\tilde{a})$ 
(see the argument given in \cite{GL}). 
Note that $2\tilde a^2 = 2\left( k+
\frac{(\ell+2)^2}{\ell} \right)$, and hence 
Theorem \ref{s3_T1} implies the estimate (\ref{s4_eq1}).
\end{lproof}

\noindent
{\bf \ref{s4}.2. Example~2.} 
Let consider the open $4$-manifold $S^2 \times \R^2$ 
with tame end $(S^2 \times S^1) \times [0, \infty)$ and 
the connected sum $X = \C\P^2 \sharp (S^2 \times \R^2)$. 
We first note that 
$$
H^2(X; \Z) \cap \Im \left[H_c^2(X; \R) \to H^2(X; \R)\right] 
\cong \ H^2(\C\P^2; \Z) \cong \Z. 
$$ 
Let $h_0 = h_S + d\tau^2$  
denote a product metric on $S^2 \times S^1$, 
where $h_S$ denotes the standard metric on $S^2$ of constant curvature one. 
Similarly to Theorem \ref{s4_T1}, 
we also obtain the following. 
\begin{Proposition}\label{s4_T2} 
With the above understood, 
$$
\begin{array}{c}
0< Y^{h\hbox{-}\cyl}(X) \leq 
\min\{12\sqrt{2}\pi,\ Y^{\cyl}_{[h+dt^2]}((S^2 \times S^1) \times \R)\} 
\end{array}
$$
for any metric $h \in \Riem^*(S^2 \times S^1)$ homotopic to $h_0$ in 
$\Riem^*(S^2 \times S^1)$. 
\end{Proposition} 
\begin{Remark} 
{\rm It is not clear for us whether a metric 
$h \in \Riem^*(S^2 \times S^1)$ 
homotopic to $h_0$ in $\Riem^*(S^2 \times S^1)$ 
satisfies the inequality 
$$ 
Y^{\cyl}_{[h+dt^2]}((S^2 \times S^1) \times \R) > 
12\sqrt{2}\pi\ (\ = Y(\C\P^2)\ ). 
$$ 
Let $h_r = h_S + d\tau^2$ denote the standard product metric on 
$S^2 \times S^1(r)$, where $S^1(r)$ stands for the circle of radius $r > 0$. 
Then we note that (cf.~\cite{P}) 
$$ 
0 < Y_{[h_S + g_0]}(S^2 \times \R^2) = 
\lim_{r \nearrow \infty}Y^{\cyl}_{[h_r + dt^2]}(S^2 \times S^1 \times \R), 
$$ 
where $g_0$ denotes the Euclidean metric on $\R^2$ and 
$$ 
Y_{[h_S + g_0]}(S^2 \times \R^2) := 
\inf_{u \in C^{\infty}_c(S^2 \times \R^2)}Q_{h_S + g_0}(u) 
$$ 
(cf.~\cite{SY}). 
If $Y_{[h_S + g_0]}(S^2 \times \R^2) > 12\sqrt{2}\pi$, 
then the estimate 
$Y^{h_r\hbox{-}\cyl}(X) \leq 12\sqrt{2}\pi$ 
(coming from Theorem \ref{s3_T1}) is effective for sufficiently large $r > 0$. 
} 
\end{Remark} 

\begin{small}

\end{small}

\vspace{20mm}

\begin{small}
{\sf
\noindent
Kazuo Akutagawa, Shizuoka University, Shizuoka, Japan
\\
e-mail: smkacta@ipc.shizuoka.ac.jp
\\
\\
Boris Botvinnik, University of Oregon, Eugene, USA
\\
e-mail: 
botvinn@math.uoregon.edu
}
\end{small}
\end{document}